%% file: cusp_preambel.tex
\documentclass{kluwer}
\usepackage[latin1]{inputenc}%
\usepackage{a4}
\usepackage{amssymb,latexsym}%
\usepackage{theorem}

\input makros.tex%
\begingroup%
\theoremstyle{plain}%
\theorembodyfont{\slshape}%
\newtheorem{satz}{Proposition}[section]%
\newtheorem{corollar}[satz]{Corollary}%
\newtheorem{lemma}[satz]{Lemma}%
\newtheorem{bemerkung}[satz]{Remark}%
\newtheorem{beobachtung}[satz]{Observation}%
\newtheorem{definition}[satz]{Definition}%
\newtheorem{zitat}[satz]{Citation}%
\endgroup%

\begingroup%
\theoremstyle{plain}
\theorembodyfont{\itshape}
\newtheorem{theorem}[satz]{Theorem}%
\endgroup%

\begingroup%
\theoremstyle{break}\theorembodyfont{\slshape}%
\endgroup%
\begingroup%
\theoremstyle{break}\theorembodyfont{\sffamily}
\endgroup%

\newcommand{\kurz}[1]{\weg{#1}}%
\newcommand{\lang}{}%
\let\variable=\epsilon\let\epsilon=\varepsilon\let\varepsilon=\variable

\long\def\weg#1{}%
\newcommand{\name}[1]{\label{#1}\textbf{Label:#1}\quad}%
\let\name\label 
\let\setminus=\smallsetminus%
\let\variable=\epsilon\let\epsilon=\varepsilon\let\varepsilon=\variable 
\begin{document}
\include{intro}%
\include{setting}%
\include{main}
\include{unipotents}
\include{hangout}%
\include{bibliography}
\end{document}

%% file: makros.tex
\def\proof{\emph{Proof: }}
\def\qed{\hspace*{\fill}$\square$}

\newcommand{\ol}{\overline}%
\newcommand{\wt}{\widetilde}%
\newcommand{\wh}{\widehat}%
\newcommand{\normaleq}{{\,\raisebox{.17ex}{\§$\vartriangleleft$\§}
    \mkern-19mu\leqslant}}%
\newcommand{\XG}{X_{\G.}}%
\def\Gn#1_#2.{\def\temp{#1} \mathbf{G}_{#2}%
      \ifx\temp\empty\else (#1) \fi}%
\def\G#1.{\Gn#1_.}%
\def\Gan#1.{\def\temp{#1} \mathbf{G}^\mathsf{an}%
           \ifx\temp\empty\else (#1) \fi}%
\def\Gis#1.{\def\temp{#1} \mathbf{G}^\mathsf{is}%
           \ifx\temp\empty\else (#1) \fi}%
\def\AdGan#1.{\def\temp{#1} \Ad\mathbf{G}^\mathsf{an}%
           \ifx\temp\empty\else (#1) \fi}%
\def\AdGis#1.{\def\temp{#1} \Ad\mathbf{AdG}^\mathsf{is}%
           \ifx\temp\empty\else (#1) \fi}%
\def\AdG#1.{\def\temp{#1} \mathrm{Ad}\mathbf{G}%
           \ifx\temp\empty\else (#1) \fi}%
\def\wtG#1.{\def\temp{#1} \wt{\mathbf{G}}%
           \ifx\temp\empty\else (#1) \fi}%
\let\jura=\S %
\def\S#1.{\def\temp{#1} \mathbf{S}%
           \ifx\temp\empty\else (#1) \fi}%
\def\wtS#1.{\def\temp{#1} \wt{\mathbf{S}}%
           \ifx\temp\empty\else (#1) \fi}%
\def\T#1.{\def\temp{#1} \mathbf{T}%
           \ifx\temp\empty\else (#1) \fi}%
\def\M#1.{\def\temp{#1} \mathbf{M}%
           \ifx\temp\empty\else (#1) \fi}%
\def\P#1_#2.{\def\temp{#1} \mathbf{P}_{#2}%
              \ifx\temp\empty\else (#1) \fi}%
\def\olP#1_#2.{\def\temp{#1} \ol{\mathbf{P}}_{#2}%
           \ifx\temp\empty\else (#1) \fi}%
\def\Ps#1_#2^#3.{\def\temp{#1} \mathbf{P}_{#2}^{#3}%
           \ifx\temp\empty\else (#1) \fi}%
\def\U#1_#2.{\def\temp{#1} \mathbf{U}_{#2}%
           \ifx\temp\empty\else (#1) \fi}%
\def\Z#1_#2:#3.{\def\Ztemp{#3} \mathcal{Z}_{#2}(#1)%
           \ifx\Ztemp\empty\else (#3) \fi}%
\def\N#1_#2:#3.{\def\Ntemp{#3} \mathcal{N}_{#2}(#1)%
           \ifx\Ntemp\empty\else (#3) \fi}%
\def\Ru#1:#2.{\def\Rutemp{#2}  R_u(#1)%
           \ifx\Rutemp\empty\else (#2) \fi}%
\def\Trans#1#2_#3.{\mathsf{Trans}_{#3}(#1,\,#2)}%
\def\Aut#1_#2.{\mathsf{Aut}_{#2}(#1)}%
\def\Stab#1_#2.{#2_{\{#1\}}}%
\def\Fix#1_#2.{{#2_{#1}}}
\newcommand{\do}{\ding{47}}%
\newcommand{\RR}{\mathbb{R}}%
\newcommand{\NN}{\mathbb{N}}%
\newcommand{\fchar}{\mathsf{char}}%
\newcommand{\Ad}{\operatorname{Ad}}%
\newcommand{\q}[1]{{\sf q}(#1)}%
\newcommand{\inner}{\mathsf{int}}%
\newcommand{\modul}{\mathsf{mod}}%
\newcommand{\rank}{\hbox{\rm rank\,}}
\def\Romannumeral#1 {\uppercase\expandafter{\romannumeral#1}}%
\def\§{\hspace*{-\the\mathsurround}}%
\newcommand{\bbackslash}{\backslash\!\!\backslash}%
\def\Splited#1>#2#3<#4{\unitlength.65mm%
\mbox{$\displaystyle #1$\rule[-3.6mm]{0mm}{3.6mm}%
\begin{picture}(25,10)(0,5.5)
      \put(2.5,7){\vector(1,0){20}}
      \put(12,8.5){$\scriptstyle #2$}
      \put(17.5,5){\oval(10,10)[br]}
      \put(17.5,0){\vector(-3,1){15}}
      \put(12.5,2){$\scriptstyle #3$}
\end{picture}
$\displaystyle #4$}}
\def\picArrow#1{\unitlength.65mm%
\begin{picture}(25,10)(0,5.5)
      \put(2.5,7){\vector(1,0){20}}
      \put(12,8.5){$\scriptstyle #1$}
\end{picture}}
\def\SES#1#2>#3>#4#5<#6{%
\mbox{$\displaystyle 1$\picArrow{}
$\displaystyle #1$\picArrow#2\Splited#3>#4#5<{#6}
\picArrow{}$\displaystyle 1$}}
\def\picarrow#1{\unitlength.35mm%
\begin{picture}(25,10)(0,9.3)
      \put(2.5,13){\vector(1,0){20}}
      \put(12,15){$\scriptstyle #1$}
\end{picture}}
\def\ses#1#2>#3>#4#5<#6{%
\mbox{$\textstyle 1$\picarrow{}
$\textstyle #1$\picarrow#2\splited#3>#4#5<{#6}
\picarrow{}$\textstyle 1$}}
\def\splited#1>#2#3<#4{\unitlength.35mm%
\mbox{$\textstyle #1$\rule[-3.00mm]{0mm}{3.00mm}%
\begin{picture}(25,10)(0,9.3)
      \put(2.5,13){\vector(1,0){20}}
      \put(12,15){$\scriptstyle #2$}
      \put(17.5,11){\oval(10,10)[br]}
      \put(17.5,6){\vector(-3,1){15}}
      \put(17.5,0){$\scriptstyle #3$}
\end{picture}
$\textstyle #4$}}


%% file: intro.tex
\begin{opening}
\title{Cusps of lattices in rank 1 {L}ie groups over local fields}

\author{Udo \surname{Baumgartner}\email{baumgart@math.uni-frankfurt.de}}

\institute{Fachbereich Mathematik\\
           Johann Wolfgang Goethe--Universit\"at\\
           60054 Frankfurt am Main}

\runningtitle{Cusps of lattices in rank 1 groups}
\runningauthor{Baumgartner}

\begin{ao}
Kluwer Prepress Department\\
P.O. Box 990\\
3300 AZ Dordrecht\\
The Netherlands
\end{ao} 

\begin{abstract} 
Let $G$ be the group of rational points of a semisimple algebraic  
group of rank 1 over a nonarchimedean local field. 
We improve upon Lubotzky's analysis of 
graphs of groups describing the action 
of lattices in $G$ on its Bruhat--Tits tree 
assuming a condition on unipotents in $G$.  
The condition holds for all but a few 
types of rank 1 groups. 
A fairly straightforward simplification of Lubotzky's definition 
of a cusp of a lattice is the key step to our results.  
We take the opportunity to reprove Lubotzky's part in the analysis 
from this foundation. 
\end{abstract}


\classification{MSC}{22E40 (20G25)}
\end{opening}

%% file: setting.tex
\section{Introduction and basics}

When properly formulated, many results true for semisimple real Lie groups 
of non-compact type 
like $SL_2(\mathbb{R})$ continue to hold, when we switch to 
(algebraic) semisimple groups over non--discrete totally disconnected 
locally compact fields (local fields for short). 

Since the topology the field induces on the group is fairly weak, 
this is surprising at first. 
It is explained by the existence of 
a combinatorial analogue of the symmetric space of a semisimple Lie group 
of non-compact type, the Bruhat--Tits building of the group. 

We are interested in lattices (discrete subgroups of finite covolume) 
in groups of rank 1. 
For groups of rank 1 the building is a tree. 
Lubotzky proved, in \cite{l.rank1}, the existence of moduli spaces of 
lattices in rank $1$ groups over local fields by giving 
recipes for their construction  
and a classification result, which describes the quotient graph of 
groups obtained from the action of the lattice on the Bruhat--Tits 
tree. 

The latter form of classification is particularly satisfactory. 
(It gives the quotient space, and at each point the stabilizer 
of a representative of the orbit mapped to this point 
drawn from a chosen connected fundamental domain.) 
This gives a concrete (and as we will see finite) method 
to construct lattices.  
The quotient graphs of lattices in rank~1 groups look like stickman 
drawings of some hyperbolic 2--manifold. 
The behaviour of the small class of arithmetic lattices in 
$SL_2$ over a function field as illustrated by the pictures 
accompanying Theorem~9 in Chapter~II.2 in \cite{trees} 
already gives you the general picture.  
 
Lubotzky's construction recipes for lattices are 
of this graph of groups flavour. 
We draw the readers attention to the method used in Theorem~5.1 in 
loc.cit. to construct nonuniform lattices from lattices in the 
unipotent radical of a minimal parabolic subgroup and a free group 
on finitely many generators. 

Our main result, part~(4) of Theorem~\ref{X mod Gamma} allows a converse 
of this to be proved. Lubotzky also offers a converse, Theorem~7.1 in 
his paper, but it is weaker, since the link back to 
lattices in appropriate unipotent radicals of parabolics is missing. 

This article could be described as the reshaping of the classification 
part of Lubotzky's paper using a simpler, yet equivalent, definition of
a cusp.  
The cusps of a lattice are supposed to describe the behaviour 
at infinity of the quotient space. We introduce them at the 
beginning of our main section, number~\ref{main}.  
In the same section, the Structure Theorem~\ref{X mod Gamma} 
is stated and its part~(4) proved. 

There is no need to prove the parts (1) and (3) which give 
already a quite precise description of the quotient graph 
of groups. 
Since our and Lubotzky's cusp definition 
are equivalent this is immediate from Lubotzky's work. 
However, many ingredients which go into Lubotzky's proof 
turn out to be not that obvious. 
We therefore decided to provide a more detailed version of 
Lubotzky's proof in our last section. 
Most of the ingredients needed will be proved in earlier sections.

Part~(2) of the Structure Theorem is an easy, but perhaps non--obvious 
application of Citation~\ref{l-6.5} and is also deferred to the last 
section. 

Most readers are expected to skip the rest of the paper, 
with the exception of Observation~\ref{anisotropListe}, which gives 
the list of types of groups who possibly do not satisfy the condition 
under which part~(4) of the Structure Theorem is true. 
We already noted, that the arguments in the last section 
are essentially known. 
The penultimate section contains the (partly difficult) 
facts on unipotent elements used. 
Most notable is an assumption which goes into an important 
ingredient of the proof of Citation~\ref{l-6.5} which is 
central in proving Lubotzky's results. 
I am grateful to Raghunathan for providing the necessary argument; see
Section~\ref{l-cit-section}. 

This leaves the rest of this section to discuss: 
We list the basic facts on   
the action of rank 1 groups on their Bruhat--Tits trees and 
of automorphism groups of trees we will need. 
In summary, the classification problem 
for lattices in rank 1 groups can be treated geometrically. 
This is only true up to a compact kernel. 
The reduction steps discussed in Section~\ref{reduction} 
take care of that problem. 

Our notational conventions are standard except that 
we write $\Fix A_G.$ for the fixator (pointwise stabilizer) 
of a subset with respect to the action of the group $G$ and  
$\Stab A_G.$ for its stabilizer.   

For constructions from the theory of algebraic groups 
we use the notation from \cite{EdM3.17} 
whose first chapter nicely summarizes most of the results 
from that theory we will need. 

\acknowledgements 
This paper owes its existence to encouragement by Alexander
Lubotzky 
and its most substantial contribution, Corollary~\ref{R-Corollar} and
Proposition~\ref{R-Satz}, to Raghunathan. 
I am grateful for their generous support.

\subsection{The geometric interpretation of a rank 1 group}
\label{rk1.geo-int}
This subsection lists the basic properties of the geometric 
approach to the analysis of lattices in the semisimple group $\G k.$. 
There is some overlap with the next section.

The group $\G k.$ acts by automorphisms on its Bruhat--Tits building 
$X(\G.,k)=:\XG$. 
It is a locally finite tree. 
The orbit of a single edge covers $\XG$, hence there are at most 
2 orbits of vertices hence at most two orders of ramification, 
$q_0+1$ and $q_1+1$. 
As another consequence, the image of $\G k.$ is cocompact in the 
group of all automorphisms of $\XG$ with its natural topology (see below). 

Each maximal compact subgroup is either the stabilizer of a vertex 
or the stabilizer of the midpoint of an edge of $\XG$. 
The latter possibility occurs iff $\G k.$ acts with inversion. 
This can only happen when $q_0$ equals $q_1$, but if this 
is the case, it will happen if $\G.$ is adjoint. 
Stabilizers of points are open. 
Each bounded subgroup of $\G k.$ fixes a point. 
Therefore the kernel of the action of $\G k.$ on $\XG$ is compact.  
The natural topology on $\Aut\XG_.$ 
makes the image of $\G k.$ in $\Aut\XG_.$ into a closed subgroup. 

From the above we derive that the image of a (uniform) lattice in 
$\G k.$ will be a (uniform) lattice in the automorphism group. 

There is a combinatorial way to compute the covolume of a lattice 
$\Gamma$ in $G:=\G k.$ (or $\Aut\XG_.$): 
The group $G$ acts without inversion on the barycentric subdivision 
$X'$ of $\XG$. The intersection of the $\Gamma$--stabilizers of 
two neighboring vertices $y_0$, $y_1$ of $X'$ thus is the 
stabilizer of the connecting edge, $e_0$ say. 
Put $U:=G_{e_0}=G_{y_0}\cap G_{y_1}$ and normalize the Haar measure 
$\mu$ on $G$ to give volume $1$ to $U$ if $\Gamma$ does not invert 
edges of $X$ and $\frac{1}{2}$ if it does. 
Then the formula for the volume of the quotient of 
$G$ modulo $\Gamma$ on page~84 in \cite{trees} gives: 
$$
\mu(G/\Gamma)=\sum_{e\in E(Y)} \frac{1}{|\Gamma_e|}\;,
$$
where $Y$ is a fundamental $\Gamma$--transversal in $X'$,  
(cf. \cite[I.(2.6)]{graphs} for notation)
and $E(Y)$ is its set of edges.

\subsection{Groups of automorphisms of trees}\label{tree-section} 

Let $X$ be a locally finite tree. 
By making each edge isometric to the unit interval 
we turn $X$ into a locally compact metric space $(X,d)$. 
Spheres and balls in $X$ with center $x$ and radius $r$ will 
we denoted $S_x(r)$ and $B_x(r)$ respectively. 
By a \emph{line (a ray)} in $X$ we mean a subspace isometric to the real line
(the set of positive reals). 
We will only use rays emanating from a vertex. 
This allows us to identify rays and lines with the sequence of vertices 
lying on them. 

Each automorphism $\alpha$ of $X$ either has a fixed point (not necessary 
a vertex) or a unique invariant line, called its axis, 
on which it induces a translation of nonzero amplitude $l(\alpha)$. 
It is called \emph{elliptic} or \emph{hyperbolic} accordingly.   

The ends of the topological space $X$ can be described purely 
combinatorially, using rays as follows: 
An \emph{end} is an equivalence class of rays. Two rays define the 
same end iff their intersection is a ray. 
For any vertex of $X$ and every end $\epsilon$ there is a unique ray, 
written $[x,\epsilon)$, emanating at $x$ and representing $\epsilon$. 
The unique point on $[x,\epsilon)$ at distance $1$ to $x$ will be written 
$\varphi_\epsilon(x)$. This map, which moves vertices closer to 
the end $\epsilon$ is useful in describing other constructions. 

Given two distinct ends, $\epsilon$ and $\epsilon'$, 
there is a unique line, called 
\emph{the line joining $\epsilon$ and $\epsilon'$} 
and written $[\epsilon,\epsilon']$, 
which contains rays representing $\epsilon$ and  $\epsilon'$ respectively. 
The points on this line are the fixed points of 
$\varphi_\epsilon^k \circ \varphi_{\epsilon'}^k$ for arbitrary $k\ge 1$. 

The ends $\epsilon$ and $\epsilon'$ of $X$ are also ends of 
the line $[\epsilon,\epsilon']$, and will be called 
\emph{the end points of the line}. 
The only ends fixed by a hyperbolic element are the end points of its axis. 

We make implicit use of the existence of a compact topology on 
the union of $X$ and its set of ends in Proposition~\ref{good/bad/ugly_geom}  
\footnote{Aside: In case $X$ is a higher rank building we have to use the 
cone topology on the geodesic compactification instead.}.   
Lets therefore describe the \emph{end topology} by describing a set of 
basic neighborhoods for any given end $\epsilon$:
Remove an edge from $X$. 
Only one of the connected components left will contain a ray 
representing $\epsilon$.  
We choose this component to be an open neighborhood of $\epsilon$.
%

Each automorphism $\alpha$ of $X$ induces a homeomorphism of its 
end compactification, and we have 
$\alpha\circ\varphi_\epsilon = \varphi_{\alpha(\epsilon)}\circ\alpha$. 

Given vertices $x$ and $y$ and an end $\epsilon$ we put 
$$
(x,y)_\epsilon:=\lim_{k\to\infty}
\left(d(x,\varphi_\epsilon^k(x))-d(y,\varphi_\epsilon^k(x))\right)\;.
$$ 
Take an arbitrary point $x$ from $X$ and define 
the \emph{oriented length function} $l_\epsilon$ 
\emph{with respect to an end $\epsilon$} by 
$l_\epsilon(\alpha):=(x,\alpha(x))_\epsilon$. 
Restricted to the stabilizer of $\epsilon$ it is 
a homomorphism to $\mathbb{Z}$ with kernel the elliptic elements. 
If $h$ is hyperbolic, and $\epsilon$, $\epsilon'$ are the end 
points of its axis, then $l_\epsilon(h)= -l_{\epsilon'}(h)$ and 
their absolute value equals $l(h)$. 
We call the end whose oriented length function with respect to $\epsilon$ 
is positive (negative) at $h$ \emph{attracting (repelling) for $h$}. 

Let $\epsilon$ be an end, $x$ a vertex and $0<\lambda\le 1$ a real number. 
Denote the unique point at distance $t$ from $x$ on $[x,\epsilon)$ by $x(t)$. 
The set 
$B_\epsilon(x;\lambda):=\bigcup_{t\ge 0} B_{x(t)}(\lambda t)$ 
is called a \emph{horoellipse} centered at $\epsilon$ with radius vertex $x$ 
and eccentricity $\lambda$. 
A horoellipse with eccentricity $1$ is written $B_\epsilon(x)$ and 
called a \emph{horoball}. 
The boundary of the horoball $B_\epsilon(x)$ is the set 
$$
S_\epsilon(x)=\{y\in X: \exists k\in \mathbb{N}\; 
                 \varphi_\epsilon^k(y)=\varphi_\epsilon^k(x)\}=
               \{y\in X: (x,y)_\epsilon=0\}\;.
$$  
It is called the \emph{horosphere} with center $\epsilon$ 
and radius vertex $x$. 
For an automorphism $\alpha$ fixing $\epsilon$ the general formula 
$\alpha(S_\epsilon(x))=S_{\alpha(\epsilon)}(\alpha(x))$
becomes $\alpha(S_\epsilon(x))=S_{\epsilon}(\alpha(x))$, 
hence $\alpha$ leaves some (every) horosphere invariant iff 
$l_\epsilon(\alpha)=0$, i.e., iff $\alpha$ is elliptic. 

The group $\Aut X_.$ is a topological group in the compact--open
topology. A base of neighborhoods of the identity is thus given by the 
sets of automorphisms fixing finitely many vertices. 
Stabilizers of vertices are seen to be pro-finite, 
hence $\Aut X_.$ in this topology is locally compact and 
totally disconnected. 
A subgroup $G$ of $\Aut X_.$ is closed iff the 
$G$--stabilizer of some (, equivalently any) vertex is 
closed (equivalently compact). 
A subgroup $\Gamma$ of $\Aut X_.$ is discrete iff the 
$\Gamma$--stabilizer of some (, equivalently any) vertex is 
finite. 

If a locally compact $\sigma$--compact group $G$ acts on $X$ with 
compact open vertex stabilizers, then the induced subgroup $\ol{G}$ 
of automorphisms of $X$ is closed, and the natural homomorphism 
$G\to \ol{G}$ is open with compact kernel. 

If $\Aut X_.$ acts with finite fundamental domain, then a subgroup 
$G$ of $\Aut X_.$ is cocompact iff $G$ acts with a finite fundamental 
domain. 

\subsection{Reduction steps}\label{reduction}

To reduce the case of a general connected semisimple $k$--group of
$k$--rank 1 to the absolutely 
simple case, we go through the following steps: 
\begin{enumerate}
\item Replace $\G.$ by its adjoint group $\AdG.$.
\item Decompose $\G.=\AdG.$ into its $k$--simple factors. 
      Since $\G.$ is supposed to have rank $1$, all but one, 
      say $\Gis.$, are anisotropic. Replace $\G.$ by $\Gis.$. 
      At this stage, $\G.$ is adjoint and $k$--simple.   
\item Choose a finite separable extension $K|k$ and a 
      $K$--group $\G.'$ such that $\G.'$ is absolutely simple 
      and the restriction of scalars $\mathbf{R}_{K|k}(\G.')$ 
      equals $\G.$. Replace $k$ by $K$ and $\G.$ by $\G.'$. 
      After this final replacement, the group will be absolutely simple. 
\end{enumerate}
(It looks like the last step should be unnecessary for groups 
of rank $1$.) 

What happens to the groups of rational points 
can be seen from the following list: 
\begin{enumerate}
\item The map $\hbox{Ad}_k$ induced by $\hbox{Ad}$ 
      on $k$--points has finite kernel, equal to 
      $\Z{\G.}_:k.=\Z{\G k.}_:.$ and compact (abelian) cokernel. 
\item Since the group of rational points of a reductive group 
      over a local field is compact if (and only if) the group 
      is anisotropic over the field, the projection onto 
      $\Gis k.$ has compact kernel, equal to $\Gan k.$, 
      where $\Gan.$ is the (in general almost) direct product 
      of the $k$--simple $k$--anisotropic factors. 
      Besides, it is continuous and open. 
\item Note, that $K$ is also local. The 
      restriction of skalars functor induces a topological 
      isomorphism between $\G.'(K)$ and $\G k.$ according to 
      \cite[I.1.7]{EdM3.17}. 
\end{enumerate} 

We list the changes to the Bruhat--Tits buildings next. 
We remind the reader that the Bruhat--Tits building is obtained 
from a natural valuation of its canonical root datum. 
\begin{enumerate}
\item For any central isogeny of reductive groups, the map induced 
      on the rational points induces an isomorphism of canonical 
      root data (attached to a maximal split torus and its image), 
      which leads to an isomorphism of Bruhat--Tits buildings.  
      The isomorphism is equivariant when the covering group is 
      made to act on the building of the covered group via the isogeny. 
\item The Bruhat--Tits building of an (almost) direct product 
      is the product of the Bruhat--Tits buildings of the factors. 
      The action is the product action. 
      The Bruhat--Tits building of a $k$--anisotropic group is 
      a single point (and vice versa). 
\item The restriction of skalars functor induces an equivariant 
      isomorphism of corresponding root data, hence of 
      Bruhat--Tits buildings. 
\end{enumerate}

The corresponding results for the spherical 
buildings over the fields involved hold as well. 
We remind the reader that the natural apartment system of 
the Bruhat--Tits building of a reductive group over a 
complete field is complete. 
As a result, the building at infinity (with its simpicial 
structure) is isomorphic to the spherical building of 
the group over the field in question. 

One can compute the kernel of the action of a connected 
semisimple $k$--group $\G k.$ on its Bruhat--Tits building explicitly. 
It is equal to $\Z{\G.}_:k.\cdot\Gan k.$. 
The kernel will therefore be trivial iff the group is adjoint 
and has no $k$--anisotropic factors. 

We summarize: 

Starting with a connected semisimple $k$--rank~$1$ group $\G.$, 
we may find a finite separable field extension $K|k$, an 
absolutely simple $K$--rank~$1$ group $\G.'$ and a canonical 
homomorphism $\pi\colon\G k.\to \G.'(K)$ with compact kernel  
(equal to $\Z{\G.}_:k.\cdot\Gan k.$) and cokernel such that 
$\G k.$ and $\G.'(K)$ have isomorphic Bruhat--Tits buildings,  
with $\G k.$ acting via $\pi$ and $\G.'(K)$ embedding  
into $\Aut\XG_.$ as a closed cocompact subgroup.   

\subsection{Ends as parabolic subgroups}
\label{U injects}\label{U HoroOp}

We finally add some results linking algebraic information on 
$\G.$ with geometric information on its Bruhat--Tits tree. 

The ends of $\XG$ can be reinterpreted as the points of its 
building at infinity. 
Since the field $k$ is complete, the natural apartment system 
of $\XG$ is complete, hence the building at infinity coincides 
(as a set) with the spherical building of $\G.$ over $k$.  
As a consequence, every end of $\XG$ \emph{is} a proper $k$--parabolic 
subgroup $\P_.$, whose stabilizer in $\G k.$ is $\P k_.$. 
In addition, the group of rational points of the unipotent radical 
of $\P_.$, $\U k_.$ say, acts simply transitively on the set of 
parabolics opposite to $\P_.$. In other words, it acts simply 
transitively on the ends different from $\P_.$.   

The maximal $k$--split torus $\S.$ whose affine apartment is 
the line $[\P_.,\P_.']$ is determined by 
$\P k_.\cap\P_.'(k)=\Z{\S.}_{\G.}:k.$. 

Let $\S.$ be a maximal $k$--split torus. 
The set of elliptic elements $\mathsf{Z}$ in $\Z{\S.}_{\G.}:k.$ 
is a compact open normal subgroup. It fixes the affine apartment 
corresponding to $\S.$ pointwise. 
It contains every other compact subgroup of $\Z{\S.}_{\G.}:k.$, 
hence is its unique maximal compact subgroup.  
The quotient $\Z{\S.}_{\G.}:k./\mathsf{Z}$ is the translation 
lattice of the affine Weyl group of $\G.$ with respect to $\S.$. 

Let $\P_.$ be a minimal $k$--parabolic subgroup of $\G.$ 
(an end of $\XG$). 
Then for any hyperbolic element $h$ in $\P k_.$ of minimal 
translation length, we have 
$$
\P k_.=(\mathsf{Z}\rtimes h)\ltimes \Ru\P_.:k.\;.
$$
where $\mathsf{Z}$ is the maximal compact subgroup of 
$\P k_.\cap\P_.'(k)$ and $\P_.'$ is the end point of the 
axis of $h$ different from $\P_.$. 
Instead of $h$ we may choose $\P_.'$ freely among the ends of 
$\XG$ different from $\P_.$. 

The group of elliptic elements in $\P k_.$ is 
$\mathsf{Z}\ltimes \Ru\P_.:k.$.

The point of intersection of a horosphere centered at $\epsilon$ 
and a line $[\epsilon',\epsilon]$ with $\epsilon'\neq\epsilon$ is unique. 
Since $\XG$ has no leaves, any point on a horosphere can be 
represented in this way, giving rise to a surjection from 
the set of ends different from $\epsilon$ to the points on a 
particular horosphere centered at $\epsilon$. 
This map is compatible with the action of the group of elliptic elements 
in the stabilizer of $\epsilon$ on both sets. 
As a consequence any group acting transitively on the ends 
different from $\epsilon$ and by elliptic transformations 
will act transitively on each of the horospheres centered at $\epsilon$. 

It can be shown, that unipotent elements always act by elliptic 
transformations (see the remarks following Theorem \ref{all good}). 
As a consequence, the group $\U k_.$ introduced above acts 
transitively on each horosphere around $\P_.\supseteq \U k_.$. 

At times we make use of the fact, that the group $\U k_.$ maps 
injectively into the group of automorphisms of $\XG$. 
This follows from two facts:\\
First, the map $\hbox{Ad}$ induces an isomorphism between the unipotent 
radical of a $k$--parabolic and the unipotent radical of its image 
(use \cite[2.15, 2.20 and 2.26]{ihes27+}). 
Second, the kernel of the action of $\G k.$ on its Bruhat--Tits building 
is $\Z{\G.}_:k.\cdot\Gan k.$ (as mentioned in the last section). 


%% file: main.tex
\section{Structure of the quotient graph of groups}\label{main}

In this section we analyze the quotient 
graphs of groups for lattices in semisimple groups of rank 1. 
It turns out that an infinite quotient graph can be explained 
in terms of the cusps of the lattice. 
This is a notion surprisingly similar to the one introduced in 
hyperbolic geometry.  
As in this classical example the notion interweaves 
geometric and group theoretic aspects.

Just as in the case of a lattice $\Gamma$ in $G:=SL_2(\mathbb{R})$, we 
expect the notion of a cusp to capture both geometrical 
properties of the quotient of the homogenous space of $G$ modulo $\Gamma$ 
and group theoretic properties of the $\Gamma$--stabilizers of these 
cusps. 

The geometric aspect involves points at infinity of the homogenous 
space and its $\Gamma$--quotient. 
Since over a local field the homogenous space is a tree, 
it is natural to assign this role to certain ends of the 
tree and all ends of the quotient graph respectively. 
For lattices in $SL_2(\mathbb{R})$ the group theoretic 
aspect involved is that $\Gamma$--stabilizers of points 
at infinity whose horoball neighborhoods map to neighborhoods 
of cusps are (maximal) unipotent groups. 
(Horoballs as introduced in Section~\ref{tree-section} will 
play a prominent role in this paper as well.) 

The first idea concerning the group theoretic aspect would 
therefore be maximal unipotent subgroups of $\Gamma$. 
Unfortunately, from a geometric point of view, unipotents in semisimple
groups over local fields do not necessarily behave as one expects 
extrapolating from the $SL_2(\mathbb{R})$. 
Unipotent elements can be grouped into three classes defined
algebraically (see Section~\ref{unipotent-section}) which can be
distinguished by their action on the Bruhat--Tits building as
illustrated by Proposition~\ref{good/bad/ugly_geom}.  
This Proposition makes clear that we should restrict ourselves to
consider good unipotent elements only. 
(An element/subgroup of $\G k.$ is good iff it is contained in 
the unipotent radical of a $k$--parabolic subgroup.  
The first part of Theorem~\ref{all good} explains, why we did not
encounter that problem for groups over the reals.) 

We thus arrive at the following definition:

\begin{definition}[cusps]\label{cusps}
Let 
$\Gamma$ be a lattice in $\G k.$. 
\begin{itemize} 
\item
An end of the Bruhat--Tits tree $\XG$ of $\G k.$ is called 
$\Gamma$--cuspidal iff there is a nontrivial good unipotent 
element in $\Gamma$ fixing $\epsilon$. 
\item 
A geometric cusp of $\Gamma$ is a an end of the 
quotient of the barycentric subdivision of $\XG$ modulo $\Gamma$. 
\item 
A cusp subgroup of $\Gamma$ is a maximal nontrivial 
good unipotent subgroup.
\item 
A cusp of $\Gamma$ is a $\Gamma$--conjugacy class of 
a cusp subgroup of $\Gamma$. 
\end{itemize}
\end{definition}

There is an obvious bijection between $\Gamma$--cuspidal ends 
and cusp subgroups of $\Gamma$: 
Since we assume that $\G.$ is semisimple of rank $1$, a nontrivial 
good unipotent is contained in a unique proper $k$--parabolic, 
hence fixes a unique end. 
Therefore every cusp subgroup of $\Gamma$ fixes a unique end 
as well, which is $\Gamma$--cuspidal by definition. 
Conversely any $\Gamma$--cuspidal end (, i.e., proper $k$--parabolic) $\P_.$ 
defines the cusp subgroup $\Gamma\cap\U k_.$, where 
$\U_.$ is the unipotent radical of $\P_.$. 
Part (2) of Theorem~\ref{X mod Gamma} will show, that this bijection 
can be pushed down to the quotient level. 

We next note that our cuspidal ends are 
what Lubotzky calls the cusps of the lattice in \cite{l.rank1}. 
Therefore we may use his results whenever convenient.

\begin{bemerkung}\label{lub-cusps}
Let 
$\Gamma$ be a lattice in $\G k.$. 
The $\Gamma$--cuspidal ends of $\XG$ are precisely the cusps of $\Gamma$ 
in the sense of Lubotzky (\cite[Definition 6.4]{l.rank1} --- see below). 
\end{bemerkung}
\proof
Recall that Lubotzky calls an end $\epsilon$ of $\XG$ a cusp of 
$\Gamma$ if and only if there exists a vertex $x$, such that writing 
$x$ as $g.x_j;\ j=0,1$ (where $x_0$, $x_1$ are representatives) 
as we may, we find a nontrivial good unipotent 
in the group $\Gamma\cap gN_1g^{-1}$ fixing $\epsilon$. 
(The definition of $N_1$ can be found just before Citation~\ref{R-3.15}; 
it depends on $j$, a natural number $n$, and a lattice $L$ in the 
Lie algebra of $\G.$. The dependence on $j$ is discussed at the beginning of 
Section~\ref{l-cit-section}. The existence of an appropriate $L$ is proved 
in the same section. The parameter $n$ must be chosen large enough to 
guarentee a geometric property $(\bigstar)$ found at the beginning of 
Section~\ref{details}.) 
Call such an end a Lubotzky--cusp for the sake of this proof. 

Let $\epsilon$ be a Lubotzky--cusp of $\Gamma$. 
By Lubotzky's definition there is even a nontrivial good unipotent 
in $\Gamma_\epsilon\cap gN_1g^{-1}$, so $\epsilon$ is evidently 
$\Gamma$--cuspidal. 
For the converse, let the end $\epsilon$ be cuspidal, i.e., 
assume there is a nontrivial good unipotent element $u$ in $\Gamma$ 
fixing $\epsilon$. 
Choose a sequence 
$(g_i)_{i\in\NN}$ in $\G k.$ such that $g_i^{-1}ug_i$ converges 
to $e$ if $i$ tends to infinity (use Theorem \ref{RPM}). 
In addition let $r\in\mathbb{N}$ be large enough to guarantee 
$N:=\Fix {B_x(r)}_{\G k.}.\subseteq N_1$. 

If we choose $k\in\mathbb{N}$ large enough to ensure 
$g_k^{-1}ug_k\in N$, we get 
$e\neq u\in\Gamma_\epsilon\cap g_kNg_k^{-1}\subseteq
  \Gamma_\epsilon\cap g_kN_1g_k^{-1}$.
We may then choose the vertex $x$ to be $g_k.x_j$ (where $j$ depends 
on $N_1$).\qed 

A further remark is in order here: 
The group $N_1$ is only defined for absolutely simple groups over 
fields of positive characteristic. 
For fields of characteristic $0$ all lattices in $\G k.$ are uniform 
thanks to a result of Tamagawa;\label{Tamagawa} see page 84 in \cite{trees}.  
But then a lattice $\Gamma$ can not contain any nontrivial good unipotent 
elements thanks to Corollary \ref{cocp=>no goods}. 
So there are no $\Gamma$--cuspidal ends hence no Lubotzky--cusps 
unless the field has positive characteristic. 
If the group is not absolutely simple, Lubotzkys definition has to be 
adapted by working backwards through the reduction steps 
of Section \ref{reduction} to cover the general case.

The following theorem gives the best discrete 
analogue of the classical description of the structure 
of the quotient of the upper half plane by a lattice in 
$SL_2(\mathbb{R})$ one can hope for. 
To guarantee that $\G k.$ acts without inversion,  
we work with the barycentric subdivision 
$\XG'$ of $\XG$ in the following theorem. 
The function $\q\cdot+1$ will give the order of ramification at 
a vertex or any vertex in the $\Gamma$--orbit of a vertex in the 
quotient graph as appropriate. 


\begin{theorem}[structure of $\Gamma\bbackslash \XG'$]
\label{X mod Gamma}
Let $\G.$ be a connected semisimple $k$--group of $k$--rank 1.
\begin{itemize}
\item[(1)] For any lattice $\Gamma$ in $\G k.$ the quotient graph 
           $\Gamma\backslash X_G'$ is the union of a finite 
           connected graph $E$ with finitely many simplicial rays 
           $r_i;\ 1\le i\le c$ attached to $E$ at their respective origin. 
 
           If $y$ and $y'$ are two neighboring vertices on one of these rays 
           which are sufficiently far from $E$ and with $y$ nearer to $E$ 
           as $y'$ then $\Gamma_y$ is a subgroup of $\Gamma_{y'}$ of Index 
           $\q{y'}$. 
\item[(2)] The map from the set of cusps of $\Gamma$ to the set of 
           geometric cusps of $\Gamma$ 
           induced by the map sending 
           each maximal good unipotent subgroup to the unique end 
           it fixes is bijective. 
\item[(3)] The $\Gamma$--cuspidal ends of $X_G'$ are precisely the ends 
           whose $\Gamma$--stabilizer is infinite and locally finite. 
           They are maximal infinite locally finite subgroups of $\Gamma$.  
           Every infinite and locally finite subgroup of 
           $\Gamma$ fixes a unique end. 
\item[(4)] Suppose that all bad unipotent elements of $\G k.$ contained 
           in $\Gamma$ are anisotropic. 
           Then each cusp subgroup of $\Gamma$ is a cocompact 
           lattice in the group of elliptic elements fixing the 
           $\Gamma$--cuspidal end $\epsilon$ fixed by the cusp subgroup. 
           In other words, the cusp subgroup is of finite index in 
           $\Gamma_\epsilon$. 
 \end{itemize}
\end{theorem}

Lubotzky (using his notion of cusps)  
proves parts (1) and (3) and states part (2).  
We will therefore confine ourselves here to report the proof of (4). 
A complete proof of the remaining parts  
can be found in Section~\ref{details}. 
 
The proof of (4) is an immediate Corollary to a technical result derived in the 
central steps~4.2--4.5 in Raghunathans proof of Citation \ref{R-4.1}. 
What makes it worth reporting, is that the condition we impose 
in (4) almost always holds, see Remark~\ref{anisotropListe}. 

We first make sure, that it suffices to treat the 
absolutely simple case: 
The reduction steps are explained in Section~\ref{reduction}. 
The facts listed in Proposition~\ref{g/b/u.reduction} 
make sure that the lattice obtained after going through the 
reduction steps still satisfies our precondition. 
On the other hand the last fact mentioned in Section~\ref{U injects} 
makes sure, that the conclusion will hold for the original 
lattice once it is proved for its ``reduced'' version. 
So we may indeed assume that $\G.$ is absolutely simple. 

If the field $k$ has characteristic $0$, we  
know (c.f. the remark following Remark~\ref{lub-cusps}) 
that there are then no nontrivial good unipotents 
in the lattice $\Gamma$ and therefore no cusp groups. 
It follows that the claim of (4) 
is trivially valid for fields of characteristic $0$. 
So we may also assume that the field has positive characteristic. 
This will enable us to use the results of \cite{R}.  

Now, let $V\leqslant \Gamma$ be a cusp subgroup. 
Let $\P_.$ be the unique end of $\XG$ it fixes, $\U_.$ its unipotent 
radical, $P:=\P k_.$ and $U:=\U k_.$. 
We have to show that $V$ is a cocompact lattice in the group 
of elliptic elements in $\P k_.=P$. It is evidently discrete. 

By definition $V$ is contained in $U$. 
Since $U$ is cocompact in the group of elliptic elements of $P$ 
(see Section~\ref{U injects}), 
it suffices to show that $V$ is cocompact in $U$.   
The unipotent group $\Lambda$ defined in the first paragraph on page~142 
of \cite{R} satisfies 
$$
V=U\cap\Gamma\leqslant \Lambda\leqslant P\cap \Gamma\;.
$$
Raghunathan proves that $\Lambda$ is cocompact in the group $L$ of 
rational points of its Zariski--closure, and that $L$ contains $U$.

Our hypotheses on unipotents in $\G k.$ enables us to prove that 
actually $V=\Lambda$ holds and our claim will follow. 
All elements of $\Lambda$ are unipotent, and none of them is
anisotropic, since they are all contained in $P=\P k_.$. 
By our hypotheses then, $\Lambda$ can not contain any bad unipotents,
therefore $\Lambda\subseteq U\cap\Gamma=V$. \qed

Assuming the additional assumption we made to derive 
part~(4) of Theorem~\ref{X mod Gamma} 
we can improve upon Lubotzky's converse to his receipee 
to construct nonuniform lattices as follows.

\begin{satz}\label{zalesskii}
Let $\Gamma$ be a lattice in a connected semisimple $k$--group of
$k$--rank 1. Assume that all cusp stabilizers of $\Gamma$ 
are residually finite and that all bad unipotent elements of $\Gamma$ 
are anisotropic. 

Then there is a sublattice $\Gamma^*$ of $\Gamma$ 
whose cusp stabilizers are cusp subgroups of $\Gamma^*$ 
and which is the free product of the representantives  
$\Delta^*_1,\ldots,\Delta^*_{c^*}$ of it's (algebraic) cusps and a free 
group of rank 
$\rank_\mathbb{Z}(H_1(\Gamma^*\backslash\XG'))$ generated by 
hyperbolic elements. 
\end{satz} 

\proof 
The method of proof applied to derive Lubotzky's converse, 
Theorem~7.1 in \cite{l.rank1}, can be reused. 
We supplement it with a geometric interpretation. 

Part (1) of the Structure Theorem   
implies that a lattice in the group of rational points of a semisimple 
group of rank $1$ over a local field is the fundamental group 
of a finite graph of groups, obtained by  ``contraction of cusps'': 
Along the simplicial rays $r_i$ we eventually have an increasing 
chain of vertex groups. 
We may therefore replace an appropriate tail of each of these rays 
of groups by a single point whose attached vertex group equals the direct 
limit (the fundamental group) of the tail. 

The fundamental group is unchanged by this modification. 
But the new graph of groups is finite with finite edge groups and 
residually finite vertex groups. 
This is obvious except for the vertex groups obtained by contraction 
of a tail of a geometric cusp. 
The latter vertex groups stabilize the cuspidal end covering the tail.  
They are therefore residually finite by assumption. 
(Indeed those vertex groups \emph{are} the stabilizers of that cuspidal end, 
since cusp stabilizers consist of elliptic elements, 
compare  Corollary~\ref {SpitzenStab-}.)
By a result of Bass--Serre theory (\cite[Proposition 12]{trees}), the 
group $\Gamma$ is residually finite and the topology of subgroups with
finite index in $\Gamma$ induces the topology of subgroups with finite
index on each of the vertex groups. 

By our second assumption and part~(4) of the Structure Theorem, the cusp
subgroups of $\Gamma$ 
have finite index in the cusp stabilizers. We may therefore choose a
normal subgroup $\Gamma^*$ of finite index in $\Gamma$ such that the
intersection with $\Gamma^*$  of the vertex groups obtained by
contraction consist of good unipotent elements and the intersection of
$\Gamma^*$  with the other vertex groups is trivial. 
The group $\Gamma^*$  determines a covering of the modified graph of
groups of $\Gamma$. 
Each of the vertex groups for $\Gamma^*$  above the contracted tails is
conjugate to the intersection of $\Gamma^*$  with the fundamental group
of the tail, hence consists of good unipotent elements. The other vertex
groups for $\Gamma^*$  are trivial. 
Hence $\Gamma^*$  is the free product of its nontrivial vertex groups
extended by the free group on the set of edges outside a maximal
subtree.  Since the graph of groups of  $\Gamma^*$ is finite, its sets of
vertex and edge groups are finite, and our claim is proved modulo the 
geometric interpretation. 

To arrive at it, we interpret the universal covering tree $\ol{X}_{\G.}$ of the
modified graph of groups of $\Gamma$ in terms of the original tree $\XG$.
According to Lemma~\ref{geo-cusp.contr}, if we choose the tails to be
contracted small enough, $\ol{X}_{\G.}$  can be realized by contracting
horoballs around cusps which are independent with respect to  $\Gamma$. 
These horoballs are then independent with respect to the subgroup
$\Gamma^*$ as well. 
The groups $\Gamma^*$ and $\Gamma$ are commensurable, hence have the
same set of cuspidal ends (this is obvious when using their
charactersisation in part~(3) of the Structure Theorem). 
The stabilizers of the independent horoballs (which are the vertices of
infinite ramification index in $\ol{X}_{\G.}$) coincide thus with the
stabilizers of the cuspidal end they contain. 
This shows that the nontrivial vertex groups of the graph of groups for
$\Gamma^*$ we considered in the previous section represent the
conjugacy classes of the stabilizers of the  $\Gamma^*$--cuspidal ends. 
We are left to confirm that each of the edges of the quotient graph of
groups of $\Gamma^*$ acts hy a hyperbolic transformation and that there
are exactly $\rank_\mathbb{Z}(H_1(\Gamma^*\backslash\XG'))$ many of
them. But this is obvious from the interpretation of the action of
$\Gamma^*$ on $\ol{X}_{\G.}$ in terms of the original action. \qed 

This stronger converse enters into the statement and proof of
Theorem~4.1 in \cite{profinite.normal}. Dependence however seems not to be
critical. 

\begin{bemerkung}
Our second assumption almost always holds, for 
we are going to show that in most groups of rank~1 each bad unipotent
element is anisotropic; compare Observation~\ref{anisotropListe}. 
On the other hand it is not known, whether cusp stabilizers of a lattice
must be residually finite, equivalently, whether lattices in groups of
rank~1 necessarily are residually finite. 
For the subclass of lattices whose cusp subgroups have
finite index in the corresponding cusp stabilizer the question becomes
whether a lattice in the group of rational points of the
unipotent radical of a semisimple group of rank~1 over a local field of
positive characteristic is residually finite. 
If the unipotent radical is abelian, it is the additive group of a
vector space over the field and hence residually finite.  
I conjecture that each discrete unipotent subgroup of an algebraic group
over a local field of positive characteristic is residually finite. 
(Discreteness is needed, since the unipotent radical of (any) 
minimal $k$--parabolic in the group $SU_3$ over an infinite field with
respect to the standard hermitian form is not residually finite.)  
\end{bemerkung}

%% file: unipotents.tex
\section{Unipotent elements}\label{unipotent-section}

\subsection{Good, bad and ugly}\label{good/bad-section}

Recall that an element $u$ in an (affine) algebraic group is called 
unipotent iff it is unipotent in some faithful rational linear 
representation. It will then be unipotent in every representation.  
It is important to note, that in positive characteristic $p$ unipotents are 
exactly the elements of $p$--power order.
 
If the underlying group $\G.$ is semisimple (or even reductive) and 
defined over $k$, the most obvious unipotents are those lying in the 
unipotent radical of a $k$--parabolic subgroup.  
An element or group contained in the unipotent radical of a 
$k$--parabolic subgroup is called $k$--good or 
more concisely \emph{good}. 
(Unipotent) elements which are not good will be named \emph{bad}.  
An element will be called \emph{anisotropic} or \emph{ very bad} 
iff the only $k$--parabolic subgroup containing that element 
is the whole group. 
In a semisimple group nontrivial anisotropic unipotents are necessarily bad. 
We will call bad elements which are not anisotropic \emph{ugly}. 
For the reduction steps we need the know, how these concepts 
behave under central isogeny, direct products and application 
of the restriction of scalars functor.

\begin{satz}\label{g/b/u.reduction}
Let $\pi\colon\wtG. \to \G.$ be a central isogeny between 
reductive $k$--groups, which is defined over $k$.
\begin{itemize}
\item[(i)] 
$\pi$ induces a bijection between the set of all 
good unipotent elements in $\wtG k.$ and the set of 
all good unipotent elements in $\G k.$ (\cite{unip2}, 2.2). 
\item[(ii)]
An image of an anisotropic element is anisotropic, 
and an image of a bad element is bad (obvious from (i)). 
\item[(iii)]
An element in an (almost) direct product is good/anisotropic, iff 
all its components are good/anisotropic. 
In an anisotropic group any element is anisotropic. 
\item[(iv)] 
The restriction of scalars functor induces bijections between 
the good, bad and anisotropic elements of both groups respectively 
over the respective fields. (This follows again from the properties 
stated in \cite[I.1.7]{EdM3.17}.)  
\end{itemize}
\end{satz}

It is possible to picture the different kinds of unipotent elements 
in a connected semisimple group $\G.$ over a local field geometrically. 
This characterization rests on the following well known 
result: 

\begin{theorem}\label{RPM}
Let $\G.$ be a connected semisimple group defined over a local 
field $k$. An element of $\G k.$ is a good unipotent element 
iff the closure of its $\G k.$--conjugacy class contains $e$. 
\end{theorem} 

We note the following important consequence. 
It hints at a link between good unipotent elements  
in a lattice and a noncompact fundamental 
domain 
(as made precise by the Structure Theorem in case of rank~1
groups). 

\begin{corollar}\label{cocp=>no goods}
Let $\G.$ be a connected semisimple group defined over a local 
field $k$. A uniform lattice in $\G k.$ does not contain 
any nontrivial good unipotent element. 
\end{corollar}
\proof 
Suppose that $\gamma$ is a nontrivial good unipotent element 
of the lattice $\Gamma$. 
By Theorem \ref{RPM} there is a sequence $(x_n)_{n\in\mathbb{N}}$
in $\G k.$ such that $x_n\gamma x_n^{-1}$ converges to $e$ as 
$n$ converges to infinity. 
Let $(\gamma_n)_{n\in\mathbb{N}}$ be the sequence constantly $\gamma$. 
Put $G:=\G k.$. Then these objects satisfy the conditions on 
the objects with the same name in Theorem~1.12 from \cite{EdM68} 
and we conclude that the images of the points $x_n$ in the space 
$G/\Gamma$ do not have a convergent subsequence. 
It follows that $\Gamma$ can not be cocompact in $\G k.$, 
a contradiction. \qed

We now derive the geometrical distinction between the different 
types of unipotent elements. 
Since we are not using this result, we only state it for 
adjoint groups without anisotropic factors.  
As is easily seen from their definition, good unipotent elements 
fix a simplex at infinity in the spherical building of $\G k.$.  
On the other hand, anisotropic elements do not and it is in 
this respect that they are evident (that is very) bad elements. 
Using the action of the group on its Bruhat--Tits building 
we derive: 

\begin{satz}\label{good/bad/ugly_geom}
Let $\G.$  an adjoint connected semisimple group defined over a local 
field $k$ without $k$--anisotropic factors. 
Then $\G k.$ may be considered a closed subgroup of the groups of 
automorphisms of its Bruhat--Tits building (use Section~\ref{reduction}) 
and we have: 
\begin{itemize}
\item
An element of $\G k.$ is \textsf{good} iff its fixed point set contains 
balls of arbitrary large diameter. 
\item
An element of $\G k.$ is \textsf{anisotropic} iff its fixed point set 
is bounded. (As will be obvious from the first part of 
the next Theorem it will be unipotent as well iff it has order 
a power of the characteristic of the field.)  
\item 
An element of $\G k.$ is \textsf{ugly} iff its fixed point set 
is unbounded but it does not contain balls of arbitrary large 
diameter. It then obviously fixes points at infinity as well. 
(The remark in parentheses of the previous item applies.)  
\item
An element of $\G k.$ is \textsf{bad} iff its fixed point set 
does not contain balls of arbitrary large diameter. 
(The remark in parentheses of the previous item applies.)  
\end{itemize}
\end{satz}
\proof By Theorem \ref{RPM} an element is a good unipotent 
iff its fixed point set contains a $\G k.$--translate of any 
ball in the Bruhat--Tits building. Since $\G k.$ acts cocompactly 
this will hold already if the fixed point set contains arbitrary 
large balls. The characterization of good unipotents claimed follows 
and so does the characterization of the bad ones. 

In geometric terms, to say that an element is anisotropic means 
nothing else than that it has no fixed points in the spherical 
building of $\G.$ over $k$. Stated equivalently, it has no 
fixed points at infinity. This will be true iff the fixed point 
set is bounded. This gives the characterization of the anisotropic 
elements. Applying the characterization of the bad elements, 
the claim on the ugly elements is established. We are done.\qed 

(If the group $\G.$ is not adjoint, or has anisotropic factors, 
the above statements have to be modified. 
Still, a good unipotent has a fixed point set containing balls of 
arbitrary large diameter. 
The converse is not true, but one can always find an element in 
any given covering, which maps to the given 
element and is good.)  

We will need the following sharper statement on fixed 
point sets of good elements available in the case of relative rank $1$. 
There is a rather obvious analogue for groups of higher rank, which 
we leave to the reader to formulate. 

\begin{lemma}[Fixed point sets of good unipotent elements]\name{Fix(u)}
Let $\G.$ be a connected semisimple group of rank $1$ over the 
local field $k$. 
The set of fixed points of a $k$--good unipotent element $u\in\G k.$ 
in $\XG$ contains a horoellipse $B_\epsilon(x;\frac{1}{3})$ for a 
suitable chosen vertex $x$ and an end $\epsilon$ (which is 
unique, if $u\neq e$).\\ 
In particular, if $h$ is hyperbolic with repelling fixed point 
$\epsilon$, the sequence of automorphisms of $\XG$ defined by 
$h^iuh^{-i}$ tends to $e$ as $i$ tends to infinity. 
\end{lemma}
\proof
\kurz{It is obvious, that the second assertion follows from the first.} 
\lang{We begin by proving the second assertion from the first:\\ 
Note that a horoellipse $B_\epsilon(x;\frac{1}{3})$ contains with each
interior point $y$ the whole horoellipse
$B_\epsilon(y;\frac{1}{3})$. We may therefore assume, replacing $x$
if necessary, that $x$ is a vertex on the axis of $h$. The 
automorphism of $\XG$ corresponding to $h^iuh^{-i}$ will then fix 
all points in the ball $B_x(\frac{i}{3})$. 
This shows, that the automorphisms $h^iuh^{-i}$ converge to the idendity
for $i$ to infinity as desired.
\kurz{In proving the first statement, we can assume $u\neq e$.}
It remains to prove, that $u$ fixes all points of a horoellipse. To
show this, assume $u\neq e$ and   
\kurz{L}\lang{l}et 
$\P_.$ be the unique minimal parabolic $k$\§--subgroup in $\G.$
containing $u$. $\P_.$ will then be the unique end fixed by $u$. 
Let $\U_.$ denote the unipotent radical of $\P_.$, and identify $\U_.$
with a root group $\U_(a).$; $a\in\Phi(\S.,\G.)$ with respect to some
maximal $k$\§--split torus $\S.$ of $\G.$ contained in $\P_.$

One of the groups $U_{a,m}$ of the canonical filtration of $U_a:=\U
k_(a).$ 
will contain $u$. Let $\alpha_{a,m}$ be the corresponding closed
halfappartment (a ray in $\XG$). Choose $x\in \alpha_{a,m}$ and let $D$
be a vector cone (\emph{chambre vectorielle} in the terminology of
\cite{ihes41}) with $x+D\subseteq \alpha_{a,m}$. 
The group $U_{x+D}$ (notation of \cite[(7.1.1)]{ihes41}) contains
$U_{a,m}$, therefore $u$. 
According to \cite[7.4.33]{ihes41} $U_{x+D}$ will then fix the
horoellipse $B_\epsilon(x;\frac{1}{3})$ pointwise (as noted in
\cite[E8]{ihes60} the scalar $\frac{1}{2}$ must be changed to
$\frac{1}{3}$). \qed

As an immediate Corollary we obtain: 

\begin{corollar}\label{SpitzenStab-}
Let $\G.$ be a connected semisimple group of rank $1$ over the local field $k$.
Let $\Gamma$ be a lattice in $\G k.$. If $\epsilon$ is a 
$\Gamma$--cuspidal end, then all elements of $\Gamma_\epsilon$ are elliptic. 
\end{corollar}
\proof 
By definition of $\Gamma$--cuspidal, $\Gamma_\epsilon$ contains a 
nontrivial good unipotent element, $u$ say. 
Suppose that $\Gamma_\epsilon$ contains a hyperbolic 
element $h$ as well. 
We may suppose then, that $\epsilon$ is a repelling fixed 
point for $h$ and apply Lemma \ref{Fix(u)}. From the second statement 
listed there we conclude, that the automorphisms corresponding to 
$h^iuh^{-i}$ converge to $e$ for $i$ to infinity. Since the image of 
$\Gamma$ in the group of automorphisms of $\XG$ is discrete, 
this implies that $h^iuh^{-i}$ acts as the
idendity for large $i$, hence $u$ does, which is impossible, since it is 
known to fix exactly one end. \qed

The following result tells us that often all unipotent 
subgroups are good. However this is not always the case, 
see \cite[3.5]{unip2} 
(note that this example works with any local field of the correct 
characteristic!). 

\begin{theorem}\label{all good}
~\\
Let $\G.$ be a reductive group defined over a perfect field 
$k$ (e.g. one of characteristic $0$). 
Then every unipotent subgroup is good [Corollaire 3.7 in \cite{unip1}]. 

Let $\G.$ be a simply connected semisimple group defined over a field 
$k$ of characteristic $p$. If $|k:k^p|\le p$, then every unipotent 
subgroup of $\G k.$ is good [Theorem 1 in \cite{el.u(GAss.sc(p))}].  
\end{theorem} 

Thanks to Corollary~1 to Propostion~4 in Chapter~I, \jura\ ~4 in 
\cite{GMW144} all local fields satisfy one of the 
conditions on the field $k$ in the above Theorem. 

Note that this implies that unipotent elements always act by elliptic
transformations, i.e., they have fixed points. 
For either the field is of characteristic $0$, any unipotent is good 
and fixes arbitrary large balls or the field has prime charcteristic, 
in which case any unipotent has finite order and therefore fixed 
points. 

As a further consequence we get the following result, which acquires 
central importance in the article \cite{R}. 
(The reference to \cite{unip1} has to be replaced by a reference to
Theorem~\ref{all good}.)

\begin{lemma}[Lemma 3.4 in \cite{R}]
If $\G.$ is a connected semisimple group over a local field $k$, then 
every unipotent element in $[\G k.,\G k.]$ is good. 
\end{lemma}

We finally use classification to produce a fairly short list of 
groups of relative rank $1$ which may contain bad unipotent 
elements which are not anisotropic. (Along similar lines one can 
compile a slightly longer list of groups of relative rank 
$1$ which may contain bad unipotents at all.)  

Suppose that $u$ is not anisotropic, hence contained in a proper 
$k$--parabolic subgroup $\P_.$. Since we assume that the rank is $1$, 
it must be minimal, in particular minimal with the property  $u\in \P_.$. 
Proposition 3.2(B) from \cite{unip2} prodcuces an 
anisotropic $u'$ in a $k$--Levi-factor $\mathbf{L}$ of $\P_.$. 
Since $\P_.$ is a proper parabolic, $\mathbf{L}$ is a reductive group 
of strictly smaller $k$--rank than $\G.$. In our case this means 
$k\hbox{--rank}(\mathbf{L})=0$ and $\mathbf{L}$ is the so called 
anisotropic kernel of $\G.$ whose type (necessarily \textsf{A}) is 
readily computed by removing the distinguished vertices form the index 
of $\G.$. Thanks to Proposition 3.2(A) of loc. cit. $u$ will be good in 
$\G.$ iff $u'$ is good in $\mathbf{L}$. 

By Theorem 2.3 of loc. cit. $u'$ will be very good hence good in
$\mathbf{L}$ if the type of this group has no connected component of the
form $\mathsf{A}_{kp-1}$, where $p$ is the characteristic of $k$. 
Using the classification table of reductive groups over local fields 
in \cite{red(local).class} together with the list of indices from 
\cite{semis.class} we see that this will always be the case, except we hit on one of the groups $\G.$ of one of the types listed below. 
Since we already ruled out the possibility of bad unipotents in 
simply connected groups we get: 

\begin{beobachtung}\label{anisotropListe}
In an absolutely almost simple $k$--group $\G.$ of $k$--rank 1 all bad 
unipotent elements will actually be anisotropic except possibly in 
one of the following cases 
($p$ denotes the characteristic of the field as usual and the type names are 
those of the first column in the classification tables of 
\cite{red(local).class}): 
\begin{itemize}
\item $p|d$, with $d\ge 2$,\quad $\G.$ not simply connected of  
  type $^dA_{2d-1}$ (with index $^1A_{2d-1,1}^{(2)}$)
\item $p=2$,  $\G.$ either 
\begin{itemize}
\item  not simply connected of type $^2A_{3}''$ 
  (index $^2A_{3,1}^{(1)}$) 
\item adjoint of type $^2C_2$ or $^2C_3$ (index 
  $C_{2,1}^{(2)}$ and $C_{3,1}^{(2)}$ respectively) 
\item not simply connected of absolute type $\mathsf{D}$, i.e., 
  $^2C$--$B_3$, ${^4D_4}$ (both with index $^2D_{4,1}^{(2)}$),  
  $^4D_5$ (with index $^1D_{5,1}^{(2)}$) or 
  $^2C$--$B_2$ (with index $^2D_{3,1}^{(2)}$).
\end{itemize}
\end{itemize}
\end{beobachtung}

Remember that our interest for this list stems from statement~(4) 
of Theorem~\ref{X mod Gamma}. 
One may not bother about the possibility that for certain groups 
over fields of characteristic $2$ cusp subgroups may not be large 
in the stabilizer of the end they fix. 
However, the type listed under the first item potentially leads 
to infinitely many exceptional characteristics. 
It would therefore be worth to explore whether ugly elements will 
really turn up for groups of this type. 
The groups in this type are strictly isogeneous to the 
group $SL_2$ of a skew field $D$ with 
degree $d^2$ (hence index $d$) over its center. 
(If $d\ge 5$, there are $\frac{1}{2}\varphi(d)$ groups in each strict 
isogeny class. 
The adjoint group of each type would be the most interesting to treat.)

\subsection{The role of good unipotents in Raghunathan's paper}
\label{R-section}

As already emphazised, we draw heavily on the results of 
Raghunathans paper \cite{R}. 
Unfortunately, we have to clarify a subtle point (see next Section) 
involving a technical construction used there. 
We therefore have to go into some detail.

We agree that in this subsection we assume 
that $\G.$ is absolutely simple and that the characteristic 
of the field $k$ is prime. 
Denote the valuation ring of $k$ by $\mathcal{O}$. 

Raghunathan works with a specific base of neighborhoods for 
the topoloy of $\G k.$, which is obtained from a full 
$\mathcal{O}$--lattice $L$ in its Lie algebra. 
The construction runs as follows: 

Fix a maximal $k$--split torus $\S.$ of $\G.$ and 
consider the decomposition of the Lie algebra of $\G.$ 
into eigenspaces \label{sfG-construction.bof}
$$
\hbox{Lie}(\G.)=\mathfrak{v}\oplus\mathfrak{z}\oplus\mathfrak{u}
$$
with respect to the action of $\S.$. 
Here $\mathfrak{v}$ and $\mathfrak{u}$ are the sum of the 
root spaces with respect to negative respectively positive 
roots of $\S.$ and $\mathfrak{z}$ is the Lie algebra of its 
centralizer. 

Denote by $\mathsf{Z}$ the unique maximal compact subgroup of 
the centralizer; for existence see Section~\ref{U injects}. 
Choose some generator $\nu$ for the spherical Weyl group in the 
normalizer of $\S.$ and put 
$\mathsf{N}:=\mathsf{Z}\cup \nu\mathsf{Z}$. 
This is a maximal compact subgroup of the normalizer of $\S.$. 

Choose $L$ now to be any full $\mathcal{O}$--lattice in 
$\hbox{Lie}(\G.)$ such that 
\begin{itemize}\label{L-conditions}
\item[(a)] $L$ is stable under the adjoint action of $\mathsf{N}$. 
\item[(b)] $L=\mathfrak{v}\cap L \oplus\mathfrak{z}\cap L 
              \oplus\mathfrak{u}\cap L$. 
\end{itemize}

In Section~\ref{l-cit-section} we will show that such an $L$ 
exists, which is even stable under a maximal compact subgroup.  

Since $\G.$ is supposed to be adjoint, we may identify it 
with its image under the adjoint representation. We put 
$$
\mathsf{G}:=\mathsf{G}(0):=\{x\in \G k.: x(L)=L\}=\G k.\cap GL(L)
$$ 
and 
$$
\mathsf{G}(i):=\{x\in \G k.: (x-\hbox{id})(L)\subseteq \pi^iL\}
$$
for $i>0$. 
It is evident, that the family of all $\mathsf{G}(i)$ for $i>0$ consists 
of compact open subgroups, which are pro--$p$, define the Hausdorff 
topology of $\G k.$ and are normal in $\mathsf{G}$.

We now give a slightly weaker version of some of the main results from
\cite{R}. 
We use the abbreviations $N_1$ and $N_2$ for the groups 
$\mathsf{G}(2n+l)$ and $\mathsf{G}(n)$ respectively. 
The parameter $n$ will be chosen later. 

\begin{zitat}[\cite{R}, Theorem 3.15]\label{R-3.15}
Let $\G.$ be an absolutely simple algebraic 
group of relative rank $1$ over the local field $k$ 
of positive characteristic. 
For any lattice $\Gamma$ in $\G k.$, there exist (integers $l$, $N_0$ and) 
a finite set $\Delta_1\subseteq \Gamma$ of nontrivial good unipotent 
elements, such that for all $g\in\G k.$ (and $n\ge N_0$) the following 
holds:\\ 
$\Gamma\cap gN_1g^{-1}\neq 1 \Rightarrow\ \hbox{there are}\ 
\delta\in \Delta_1 \ \hbox{and}\ \gamma\in\Gamma\ \hbox{such that}\    
\theta:=\gamma\delta\gamma^{-1}\in gN_2g^{-1}$.
\end{zitat}

The following result is an easy Corollary. 

\begin{zitat}[\cite{R}, Corollary 3.16]\label{R-3.16}
Let $\G.$, $k$  and $\Gamma$ be as in the last result. 
If $\P_.$ is a $k$--parabolic subgroup of $\G.$ 
which is $\Gamma$--cuspidal, then there is a $\delta$ in $\Delta_1$
such that $\P_.$ is $\Gamma$--conjugate to the unique 
$\Gamma$--cuspidal $k$--parabolic containing $\delta$.  
\end{zitat}


We may reformulate \cite[theorem 4.1]{R} into a first general 
result on the structure of cusp stabilizers: 

\begin{zitat}[\cite{R}, Theorem 4.1]\name{R-4.1}
Denote by $^\circ\P k_.$ the kernel of the modular function of 
$\Ru\P_.:k.$ restricted to conjugation by elements of $\P k_.$.\\ 
If the $k$--parabolic  $\P_.$ is $\Gamma$--cuspidal then 
$\Gamma\cap{^\circ\P k_.}$ is a cocompact lattice in ${^\circ\P k_.}$.
\end{zitat}

\subsection{Implications and reinterpretations}

Reading Corollary~\ref{R-3.16} differently, we get for free: 

\begin{satz}\label{f of cusps}
The elements of $\Delta_1$ represent all conjugacy classes of cuspidal
ends. Hence every lattice has only finitely many cusps. 
\end{satz}

To make use of Citation \ref{R-4.1}, we have to reinterpret  
$^\circ\P k_.$ geometrically. 
Write $\modul_G$ for the modular function of a locally compact 
group. We have: 

\begin{satz}\name{0P=EP}
If $\P_.$ is a minimal $k$--parabolic subgroup of $\G.$ then 
$^\circ\P k_.$ is the group of elliptic elements of $\P k_.$; 
in fact: $\modul_{\U k_.}(\inner(h))=
(q_0q_1)^{\frac{1}{2}l_{\P_.}(h)}$ for any $h\in\P k_.$.
\end{satz}
\proof 
Since the canonical map $\G k.\to \Aut{X_G}_.$ induces an
isomorphism $\U k_.$ to it's image $U$ (c.f. Section \ref{U injects}),
we may compute the module of inner conjugation by $h\in\P k_.$
`"geometrically"' inside $\Aut{X_G}_.$ : 
Let $\mu$ be a left invariant Haar measure on $U$. We first show 
that all elliptics in $\P k_.$ are contained in 
${^\circ\P k_.}$, i.e. $\modul(\inner(g)_{|\U k_.})=1$ 
for any elliptic  $g\in\P k_.$. Chose a vertex $x$ fixed by
$g$. Using the fact that $U_x$ is an open compact subgroup of $U$ we get :
  $$\modul(\inner(h)_{|U})=
  \frac{\mu\bigl(\inner(g)(U_x)\bigr)}{\mu\bigl(U_x\bigr)}  =1\ ,$$
To compute $\modul(\inner(h)_{|U})$ for hyperbolic $h\in \P k_.$ it is
sufficient to stick to the case where $\P_.$ is attracting for $h$.
Let $x$ denote a vertex on the axis of $h$. 
We have $hU_xh^{-1}=U_{h.x}\geqslant U_x$ and 
$$
\modul(\inner(h)_{|U})=|U_{h.x}:U_x|=\#\bigl(S_{h.x}(l(h))\cap
  S_{\P_.}(x)\bigr)\ .
$$ 
Only the second equality requires proof. $U_{h.x}$ acts on $S_{h.x}(l(h))\cap
S_{\P_.}(x)$, the stabilizer of $x\in S_{h.x}(l(h))$ being $U_x$. We
are therefore obliged to show that this action is transitive. 
The group $U$ acts transitively on each horosphere with center $\P_.$ 
(Section \ref{U HoroOp}), hence in particular on $S_{\P_.}(x)\supseteq
S_{h.x}(l(h))\cap S_{\P_.}(x)$. Any element $u\in U$ mapping a point
in the latter set to another one must fix $h.x$, which proves
equality. Since the cardinality of $S_{h.x}(l(h))\cap S_{\P_.}(x)$ is
easily seen to be $(q_0q_1)^{\frac{1}{2}l(h)}$ we are done. \qed

\subsection{Interpreting Raghunathan's neighborhood basis}
\label{l-cit-section}
We need to comment on one last result needed in the proof 
of parts~(1) to (3) of Theorem~\ref{X mod Gamma}. 
It asserts that 
whenever a nontrivial good unipotent element $\theta$ of a lattice 
fixes a sufficiently large ball around some point $x$, 
then the whole stabilizer of $x$ will fix the 
end fixed by $\theta$. 
(This will be combined with Citation~\ref{R-3.15} which 
states that whenever there is a nontrivial element of the 
lattice, which fixes a huge ball around a point $x$, 
then there is a nontrivial good unipotent element fixing 
a large ball around $x$.) 

In \cite{l.rank1}, where this result was proved, 
two points were overlooked. 
For one, the statement depends on the point $x$ obviously, 
but the condition imposed varies only with the $\G k.$--orbit 
of $x$. 
There is an easy fix for this problem; 
only vertices and probably midpoints of edges are of interest 
as choices for $x$ and we simply make the data dependent on 
the corresponding $\G k.$--orbit. 
The second problem is serious. 
Since it is visible only on close inspection, we 
proceed to the precise form of the statement in question. 

Let $x_0$ be a fixed vertex, and $x_1$ some point closest to but distinct 
from $x_0$, whose stabilizer is also maximal compact. 
Note that $x_1$ will be the midpoint of an edge 
incident with $x_0$ if the tree is regular, not a vertex. 
(In that case the adjoint group acts with inversion.) 

Recall Raghunathans neighborhood base introduced at the 
beginning of Section~\ref{R-section}. 
The fix to the second problem consists in proving:\\
$(*)$\qquad There exist $\mathcal{O}$--lattices $L[0]$ and $L[1]$ 
            such that the corresponding group\linebreak
\hphantom{$(*)$\qquad }
            $\mathsf{G}[0]$ respectively $\mathsf{G}[1]$ is the full 
            stabilizer of $x_0$ and $x_1$ respectively.\\
This amounts to a geometrical interpretation of Raghunathans 
neighborhood basis. 
(The derivated groups $\mathsf{G}(i)$, $N_1$ and $N_2$ introduced 
 in Section~\ref{R-section} also aquire dependence on $x_j$.  
 They should therefore be written accordingly, but we only do this 
 when dependence on $j$ is important.)

The precise version of Lubotzkys Lemma~6.5 is then 
is follows: 
(Remark: It is obvious from the proof, that 
it is unnecessary to nail down $N$ to be $N_2$.)

\begin{zitat}[\cite{l.rank1}, Lemma 6.5]\name{l-6.5}
Let $\G.$ be an absolutly simple algebraic 
group of relative rank $1$ over the local field $k$ 
of positive characteristic and let $\Gamma$ be a lattice in $\G k.$. 
Let $x=g.x_j$ with $j\in\{0,1\}$ and $g\in\G k.$. 
Let $N=N[j]$ be (one of the groups $\mathsf{G}(i)[j]$ with $i>0$, e.g.) the 
group $N_2[j]$.\\ 
If $\Gamma\cap gNg^{-1}$ contains a nontrivial good unipotent, 
$\theta$ say,  
then writing $\epsilon$ for the unique end fixed by $\theta$ we'll have 
$\Gamma_x\subseteq \Gamma_\epsilon$. 
\end{zitat}

Lets pause for a moment, to see why we need a result like $(*)$: 
The proof of Lemma~6.5 in \cite{l.rank1} is completed by the 
observation that $Q$ is normal in $\Gamma_x$. 
This follows from $N_2\normaleq R$, 
where $R$ denotes the stabilizer of a point of the appropriate 
type (red). 
The last claim is stated as a fact in \cite[6.2]{l.rank1}. 
The reference given is \cite[3.15]{R} cited herein as Citation~\ref{R-3.15}.
As the reader can check, within \cite{R} (and \cite{l.rank1}) 
$N_2=\mathsf{G}(n)[j]$ is only shown to be normal in $\mathsf{G}$. 

But we are safe, once we show that the lattice $L$ on which 
$\mathsf{G}$ depends can be chosen as to guarentee 
that $\mathsf{G}$ equals any maximal compact subgroup we prescribe. 
That's what we are going to do now.  

Set $M$ equal to the stabilizer of $x_j$  for $j=1,\,2$. 
Take any maximal split torus $\S.$ guarenteed to exist by the Corollary 
below as the split torus needed for the construction of 
$\mathsf{G}$, starting on page~\pageref{sfG-construction.bof}. 
This torus is provided by the Proposition. 
From the proof of the latter it is clear, that the affine 
apartment of any torus which qualifies will contain $x_j$. 
(Alternatively, if we assume that the residue field of $k$ 
has at least $4$ elements, we can argue as follows: 
Maximality of $M$ implies that $M$ contains the group of units 
of the torus $\S.$. But under our assumption, the only 
fixed points of the group of units of $\S.$ are the points 
of its affine apartment thanks to \cite[3.6.1]{red(local).class}.) 

Statement (B) below then clearly implies condition (b) 
for $L$ (to be found on page~\pageref{L-conditions}). 
Statement (A) translates into $M\subseteq \mathsf{G}$. 
Since $M$ is maximal, this implies $M=\mathsf{G}$. 

The group $\mathsf{N}$ which condition (a) on the same page refers to 
is a compact subgroup of the normalizer of $\S.$. 
We will show that it can be chosen to lie in $M=\mathsf{G}$. 
Since $x_j$ lies in the apartment corresponding to $\S.$, 
we know that $\mathsf{Z}$ is contained in $M$. 
Now choose $\nu$ in such a way, that $x_j$ is its only fixed 
point in the affine apartment for $\S.$. 
Then we'll have $\mathsf{N}\subseteq M$. 

We thus are reduced to showing the Proposition and its Corollary.  
We assume that $\G.$ is a connected semisimple algebraic group and 
$p\colon\widetilde{\mathbf{G}} \rightarrow \mathbf{G}$ its universal covering. 
The weight space corresponding to a root $\alpha$ will be denoted 
Lie $\mathbf{G}(\alpha)$. The set of rational points, 
$\mathrm{Lie}\mathbf{G}(\alpha)(k)$ will be denoted Lie$G(\alpha)$.  
As before we denote by $\mathfrak{z}:=\mathfrak{z}(\S.)$ the Lie Algebra of 
the centralizer of the torus $\S.$. 

\begin{corollar}\label{R-Corollar} 
Let ${\mathbf G}$ be a 
connected semisimple algebraic group defined over the local field $k$. 
For any maximal compact subgroup $M$ of $\G k.$ there is a 
full $\mathcal{O}$--lattice $L$ in $\hbox{Lie}(\G.)(k)$ and a maximal 
$k$--split torus $\mathbf{S} \subseteq \mathbf{G}$ such that 
\begin{itemize}
\item[(A)] $L$ is stable under the adjoint action of $M$. 
\item[(B)] $L= L \cap \mathfrak{z}(\widetilde{S}) \oplus 
 {\displaystyle{\coprod_{\alpha\in\Phi(\widetilde{\S.},\wtG.)}}} 
(\mathrm{Lie}G(\alpha) \cap L)$.
\end{itemize}
\end{corollar}

In the Proposition we will make use of several field extensions 
$k'$, $\wt{k}$ and $\wh{k}$. 
The corresponding rings of integers and residue fields will be
written $\mathcal{O}'$, $\wt{\mathcal{O}}$, $\wh{\mathcal{O}}$ 
and $F'$, $\wt{F}$, $\wh{F}$  
respectively. 
The Galois group of a field extension $K|k$ will be written 
$\mathcal{G}(K|k)$. 

\begin{satz}\label{R-Satz}
Let $M$ be any compact subgroup of $G= \mathbf{G}(k)$ and 
$\wh{k}|k$ a finite unramified extension.  
Then there is a maximal $k$--split torus 
$\widetilde{\mathbf{S}}\subset\wtG.$ 
and a parahoric subgroup $\widetilde{P}$ of $\widetilde{G}= 
\widetilde{\mathbf{G}}(k)$ with the following properties. 
\begin{itemize}
\item[(i)]
$p(\widetilde{\mathbf{P}}(\wh{\mathcal{O}}))$ and $M$ generate a 
compact subgroup of $\mathbf{G}(\wh{k})$.  
(Here $\widetilde{\mathbf{P}}$ is the parabolic group 
scheme over $\mathcal{O}$ defined as spec $ R_{\widetilde{P}}$ with 
$$
R_{\widetilde{P}} = \{f \in k [ \widetilde{\mathbf{G}}] \mid 
                    f(\widetilde{P})\subseteq \mathcal{O}\}
$$ 
$k [ \widetilde{\mathbf{G}}]$ being the coordinate ring of 
$\widetilde{\mathbf{G}})$.
\item[(ii)]
The maximal compact subgroup of $\widetilde{\mathbf{S}}(\wh{k})$ 
is contained in $\widetilde{P}$. 
\end{itemize}
\end{satz}

\emph{Proof (of Corollary): }
Let $\wh{k}$ be an unramified extension chosen such that for any 
$k$--split maximal torus $\widetilde{\mathbf{S}}$ of 
$\mathbf{G}$, there is an element $t\in\widetilde{\mathbf{S}}(\wh{k})$ 
such that $\alpha(t), \alpha(t)-1$ are units for every 
$\alpha\in\Phi(\widetilde{\mathbf{S}},\wtG.)$ and $\alpha(t)-\beta (t)$ 
is a unit for every pair of distinct roots $\alpha, \beta$ 
in $\Phi(\widetilde{\mathbf{S}},\wtG.)$ 
(the residue field of $\wh{k}$ needs to be sufficiently large to secure this).  
Choose $\widetilde{\mathbf{S}}$  and a parahoric subgroup $\widetilde{P}$ 
as in the proposition (for the compact group $M$).  
Since $M$ and $p(\widetilde{\mathbf{P}}(\wh{\mathcal{O}}))$ generate a compact 
subgroup $M'$ of $\mathbf{G}(\wh{k})$, 
there is $\wh{\mathcal{O}}$--lattice $L$ in 
Lie$\mathbf{G}(\wh{k})$ which is stable under $M'$.  
The element $t$ necessarily belongs to the maximal compact subgroup 
of $\widetilde{\mathbf{S}}(\wh{k})$ (all its eigenvalues in the adjoint 
representations are assumed to be units).  By our choice of $t$, 
the eigen spaces of $\hbox{Ad}(t)$ are the same as those of 
$\widetilde{\mathbf{S}}(\wh{k})$ and it follows from elementary linear 
algebra that  
$$
L=L\cap \mathfrak{z}(\widetilde{S}) \oplus 
  \coprod_{\alpha\in\Phi(\widetilde{\mathbf{S}},\wtG.)}
  (Lie G(\alpha)\cap L)\ .$$  
Hence the corollary. \qed 

\emph{Proof (of Proposition): }
To prove the proposition we need to use Bruhat Tits theory. 
Fix a $k$--split torus $\widetilde{\mathbf{S}}$ in $\widetilde{\mathbf{G}}$ 
and a maximal torus $\widetilde{\mathbf{S}}_1$ of $\widetilde{\mathbf{G}}$ 
defined over $k$ and containing $\widetilde{\mathbf{S}}$ and such that 
$\widetilde{\mathbf{S}}_1$ contains a maximal split torus over the 
maximal unramified extension $\widetilde{k}$ of $k$. 
Suppose now that $k'|k$ is an unramified extension such that the maximal 
split torus over $k'$ contained in $\mathbf{S}_1$ is also maximal split over 
$\widetilde{k}$. Let $\mathcal{B}$ be the Bruhat-Tits building associated 
to $\widetilde{\mathbf{G}}(k')$ and $\mathcal{A}\subset\mathcal{B}$ 
the $U_1$ fixed points in $\mathcal{B}$ where $U_1$ is the maximal compact 
subgroup of $\widetilde{\mathbf{S}}_1(k')$.  
Then $\mathcal{A}$ is  ``an apartment'' in $\mathcal{B}$ and is stable 
under the action of the Galois group $\mathcal{G}(k'/k)$; 
and by the fixed point theorem of Bruhat--Tits, there is a point 
$b\in\mathcal{B}$ fixed by $\mathcal{G}(k'/k)$.  
Let $\widetilde{P'}$ be the subgroup of $\widetilde{G}(k')$ that fixes 
the point $b$. Then the $\mathcal{O}$--algebra 
$R_{\widetilde{P}}= \{ f\in k [\widetilde{\mathbf{G}}] \mid 
   f(\widetilde{P'}) \subset \mathcal{O}'\}$  
defines a parahoric group scheme 
$\widetilde{\mathbf{P}}=$ spec $R_{\widetilde{P'}}$ over $\mathcal{O}$ 
and $\widetilde{\mathbf{P}}(\mathcal{O}')=\widetilde{P'}$.  
Let $F$ be the residue field of $\mathcal{O}$ and 
$\pi=\widetilde{\mathbf{P}}(\mathcal{O}')\rightarrow 
 \widetilde{\mathbf{P}}(F)$ 
the natural map.  
$\widetilde{\mathbf{P}}\otimes_\mathcal{O} F$ is a connected group scheme over 
the finite field and hence admits a Borel subgroup $B$ over $F$; 
then $\pi^{-1} (B(F))$ (resp $\pi^{-1} (B(F'): \pi$ also denote the map 
$\widetilde{\mathbf{P}}(\mathcal{O}') \rightarrow \widetilde{\mathbf{P}}(F')$, 
is an Iwahori subgroup 
$\widetilde{I}$ of $\widetilde{\mathbf{G}}(k)$ 
(resp $\widetilde{I'}$ of $\widetilde{G}(k'))$.  
If $R_{\widetilde{I}} = \{f\in k [\widetilde{\mathbf{G}}] \mid 
f(\widetilde{I}) \subset \mathcal{O}\}$, then 
spec $R_{\widetilde{I}}= \widetilde{\mathbf{I}}$ is an 
Iwahori group scheme with 
$\widetilde{{\mathbf I}} (\mathcal{O}) = \widetilde{I}$ and 
$\widetilde{\mathbf{I}}(\mathcal{O}')= \widetilde{I'}$. 
The torus $\widetilde{\mathbf{S}}$ being split has a natural definition 
over $\mathbb{Z}$ and hence over $\mathcal{O}$. 
Moreover it is easy to see that the restriction map of functions in 
$R_{\widetilde{\mathbf{P}}}$ to $\widetilde{\mathbf{S}}$ 
gives a closed immersion of this split torus over $\mathcal{O}$ in 
$\widetilde{\mathbf{P}}$.  
From the considerations it is easy to see that $B$ above can be so chosen 
that it contains the reduction modulo the maximal ideal of the split torus 
$\widetilde{\mathbf{S}}$ (over $\mathcal{O})$. 
This means that $\widetilde{I'} = \widetilde{\mathbf{I}}(\mathcal{O}')$ 
contains the maximal compact subgroup 
$(= \widetilde{\mathbf{S}}(\mathcal{O}'))$ of $\widetilde{\mathbf{S}}(k')$.  
This leads us to the conclusion that {\it any} Iwahori subgroup of 
$\widetilde{\mathbf{G}}(k')$ stable under $\mathcal{G}(k'|k)$ necessarily 
contains the maximal compact subgroup of $\widetilde{\mathbf{S}}(k')$ 
with $\widetilde{\mathbf{S}}$ a maximal $k$--split torus in 
$\widetilde{\mathbf{G}}$. 

Suppose now that $M$ is as in the Proposition. We assume, as we may, 
that $M$ is a maximal compact subgroup of $G$. 
The group $M$ as well as the group $\mathcal{G}(\wh{k}/k)$  
($\wh{k}$ as in the proposition ) act as isometries of the Bruhat Tits bundle 
$\mathcal{B}_o$ associated to $\wtG.(\wh{k})$.  
Since their actions on $\mathcal{B}_0$ commute with each other, they generate 
together a compact group of isometries of $\mathcal{B}_0$ and hence by the 
Bruhat--Tits fixed point theorem they have a common fixed point $b$.  
Such a fixed point determines a parabolic subgroup $\widetilde{P}$ over $k$ 
which necessarily contains an Iwahori subgroup $\widetilde{I}$ over $k$.   
Let $\widetilde{\mathbf{I}}$ (resp $\widetilde{\mathbf{P}})$ be the 
group scheme over $\mathcal{O}$ associated to $\widetilde{I}$ 
(resp $\widetilde{P})$.  
Then the isotropy group in $\G \wh{k}.$ of the point 
$b \in\mathcal{B}$ contains $M$
as well as $p(\widetilde{\mathbf{P}}(\wh{\mathcal{O}}))$ hence 
$p(\widetilde{\mathbf{I}}(\wh{\mathcal{O}})))$; and from the 
preceding paragraph, 
we know that $\widetilde{\mathbf{I}}(\wh{\mathcal{O}})\supset$ 
maximal compact of 
some $\widetilde{\mathbf{S}}(\wh{k})$ with $\widetilde{\mathbf{S}}$ 
maximal split torus over $k$. This proves the proposition. \qed

%% file: hangout.tex
\section{The details}\label{details}

In this last section, we supply a complete proof of 
parts~(1) to (3) of ``our'' main result, Theorem~\ref{X mod Gamma}. 
The reasoning follows closely Lubotzkys original proof. 
At some places we have to inject facts we piled up in 
earlier sections for which \cite{l.SL2} provides an 
easier proof in the special case of $\G.=\mathbf{SL_2}$.

We start with a reduction to the case of an absolutely simple 
group, using the steps listed in Section~\ref{reduction}. 
Statements (1)~and (3) will follow in general, once they are 
proven in the absolutely simple case for the simple reason 
that the restriction of the map 
$\G k.\to \G.'(K)\leqslant \Aut\XG'_.$ 
to $\Gamma$ has finite kernel. 
Claims (2)~and (4), dealing with good unipotents,  
need in addition that the group of rational 
points the unipotent radical of a $k$--parabolic 
injects into $\G.'(K)\leqslant \Aut\XG'_.$, as stated in 
Section~\ref{U injects}. 
The reader may wish to have a glance at the detailed 
description of the map from algebraic to geometric cusps 
given in the final step of the proof of Corollary~\ref{finish} 
on page~\pageref{spell gp>geo} to verify that statement for 
claim~(4). 

Further, all claims in Theorem~\ref{X mod Gamma} are trivial for 
cocompact lattices. 
By Tamagawas result, already quoted on 
page~\pageref{Tamagawa}, we may therefore assume that the 
field $k$ has positive characteristic. 
In sum, all the results of Raghunathans paper \cite{R} 
will be applicable. 
The results to follow however will usually be true in general, 
as will be obvious a posteriori. 

Recall how the points $x_0$ and $x_1$ were chosen 
(beginning of Section~\ref{l-cit-section}): 
$x_0$ is some fixed vertex, and $x_1$ closest to but distinct from $x_0$, 
with the property that it has a maximal compact stabilizer as well. 
We now need to choose the free parameters $n=n[0],\,n[1]$ with $n\ge N_0$ 
which determines the groups $N_1$ and $N_2$ (of type $[0]$ and $[1]$) 
and a radius $\rho$ large enough such that we have 
$$\label{bigstar}(\bigstar)\quad
\Fix{B_{x_0}(\rho)}_G.\cup\Fix{B_{x_1}(\rho)}_G.\subseteq N_1[0]\cap N_1[1]
\ \ \hbox{and}\ \ 
N_2[0]\cup N_2[1] \subseteq\Fix{B_{x_0}(1)}_G. \cap \Fix{B_{x_1}(1)}_G.
$$
This is clearly possible, since the family of groups $\mathsf{G}(i)$; 
$i\in \mathbb{N}$ is a system of neighborhoods of $e$ in $\G k.$ defining the 
Hausdorff--topology and since the (adjoint, absolutly almost simple) 
group $\G k.$ can be identified with a topological subgroup of $\Aut\XG_.$. 

The following fundamental Lemma describes what happens near  
a $\Gamma$--cuspidal end. 
We will make the barycentric subdivision $\XG'$ into a metric subspace of the
metric space $\XG$ by asigning length $\frac{1}{2}$ to the edges of
$\XG'$. When considered as a function on $\XG'$ $\varphi_\epsilon$ will
translate by the distance $\frac{1}{2}$ accordingly.

\begin{lemma}\label{note1}
Let $\epsilon$ be a $\Gamma$--cuspidal end. 
Then the following assertions hold :
\begin{itemize}
\item[(i)] If a vertex $y$ is chosen to lie sufficently close to $\epsilon$,
  then 
\begin{enumerate}
\item the action of $\Gamma_\epsilon$ on $S_\epsilon(y)$ and therefore
  also on $S_\epsilon(\varphi_\epsilon^k(y))$; $k\ge 0$ will be transitive.
\item For any $y$ as in 1. \quad $\#\Gamma_{\varphi_\epsilon^k(y)}\ge
  (q_0q_1)^{\lfloor \frac{k}{3}\rfloor}$.
\end{enumerate}
\item[(ii)] Given $r\in\RR$ we can find a vertex $y:=y(r)$ sufficently
  close to $\epsilon $ such that
\begin{enumerate}
\item $\Gamma_\epsilon$ contains a non trivial good unipotent $\theta_y$
  fixing $B_{y_k}(r)$  pointwise with
  $y_k:=\varphi_\epsilon^k(y)\in [y,\epsilon[$; $k\in\NN$.    
\item if we chose $r\ge\rho$, then $\Gamma_{\varphi_\epsilon^k(y)}\subseteq
  \Gamma_{\varphi_\epsilon^{k+1}(y)} \subseteq \Gamma_\epsilon$ for all
  $k\ge 0$.  
\end{enumerate}
\item[(iii)] If we chose $y\in {X_G'}^0$ to have both properties (i).1
  and (ii).1 with $r\ge\rho$, then 
\begin{enumerate}
\item any vertex $x\in B_\epsilon(y)^0$ will also have property (ii).1
  with respect to $r$, i.e. there exists a non tivial good unipotent 
  $\theta_x\in\Gamma_\epsilon$ fixing $B_x(r)$ pointwise. Furthermore
  for any such $x$: $\Gamma_x\subseteq
  \Gamma_{\varphi_\epsilon(x)}\subseteq \Gamma_\epsilon$.
\item Any vertex $x\in B_\epsilon(y)$ has the property  
      $$\forall k\ge0\quad |\Gamma_{\varphi_\epsilon^{k+1}(x)} :
      \Gamma_{\varphi_\epsilon^{k}(x)}| = \q{\varphi_\epsilon^{k+1}(x)}\
      .$$     
\end{enumerate}
\end{itemize}
\end{lemma}
\proof 
We may suppose without loss of generality that $\G.$ is absolutly almost
simple and adjoint 
(and for trivial reasons that furthermore $\fchar k\neq 0$). 
We begin by proving (i). $\Gamma_\epsilon$ acts on the horospheres
with center $\epsilon$ according to Corollary \ref{SpitzenStab-}. 
Choose a ray $(y_i)_{i\in\NN}$ representing $\epsilon$. Regarding
$\epsilon$ as a parabolic $k$\§--subgroup $\P_.$ we find from
proposition \ref{0P=EP} that $^\circ\P k_.$ also acts on the horospheres
with center $\P_.$. This action is transitive, since the same already
holds for the subgroup $\Ru\P_.:k.$. 
${^\circ\P k_.}$  is covered by the increasing sequence of compact open
subgroups $\bigl({^\circ\P k_.}_{y_i}\bigr)_{i\in\NN}$.  Any compact
subset of ${^\circ\P k_.}$ in particular a compact system of
representatives of $\Gamma\cap {^\circ\P k_.}$ in ${^\circ\P k_.}$ 
(which exists according to \ref{R-4.1}) is contained in one of these say
in ${^\circ\P k_.}_{y_{i_1}}$.\\ 
From ${^\circ\P k_.}=(\Gamma\cap {^\circ\P k_.})\cdot {^\circ\P
  k_.}_{y_{i_1}}$ we conclude that $\Gamma\cap{^\circ\P k_.}$ acts
transitively on $S_{\P_.}(y_{i_1})=S_\epsilon(y_{i_1})$. 
The same will then be true for the horospheres
$S_\epsilon(\varphi_\epsilon^k(y_{i_1}))=S_\epsilon(y_{i_1+k})$ for all
$k\ge 0$. (i).1 will then follow with $y:=y_{i_1}$. 

Choose a vertex $y$ which has property (i).1. An element of
$\Gamma_\epsilon$ permuting the points of $S_k:=S_{\epsilon}(y) \cap
S_{\varphi_{\epsilon}^k(y)}(k)$ necessarily fixes
$\varphi_{\epsilon}^k(y)$. This implies 
$$\#\Gamma_{\varphi_{\epsilon}^k(y)}\ge \# S_k \ .      $$
To prove (i).2 it therefore remains to bound the size of $S_k$ from
below by $(q_0q_1)^{\lfloor\frac{k}{3}\rfloor}$ which is easy. 

We now turn to assertions (ii) and (iii): Choose a nontrivial good 
unipotent $\theta\in\Gamma_\epsilon$. It fixes some horoellipse
$B_\epsilon(y_0;\frac{1}{3})$; c.f. Lemma \ref{Fix(u)}. 
If $y=\varphi_\epsilon^k(y_0)$ is a vertex which so close
to $\epsilon$ to assure that on the one hand $\Gamma_\epsilon$ acts
transitively on $S_\epsilon(y)$ and on the other hand $k\ge 6r$ then we
find that $\theta=:\theta_y$ fixes $\epsilon$ and all points of the ball
$B_y(r)$. All $y'\in[y,\epsilon[$ will have the same property. We have
shown (ii).1.\\  
Before turning to the second part of (ii) we consider now (iii).1: 
Let $x\in B_\epsilon(y)$ be a vertex on the same horosphere as
$y$. Choose $\gamma\in\Gamma_\epsilon$ mapping $y$ to $x$ and put
$\theta_x:=\gamma \theta_y\gamma^{-1}\in\Gamma_\epsilon$. This is
evidently a good unipotent element suitable for all vertices
$\varphi_\epsilon^k(x);\;k\ge 0$. Since the rays $[x,\epsilon[;\; x\in
S_\epsilon(y)$ cover the whole horoball $B_\epsilon(y)$ we checked all
parts of (iii).1 modulo (ii).2. 

To proof (ii).2 choose a vertex $y$ with property (ii).1. 
To handle $y_k\in[y,\epsilon[$ choose 
$g_k$ mapping either one of $x_0$, $x_1$ to $y_k$ 
arbitrarily. $\theta_y=\theta_{y_k}$ will then fix the ball
$B_{y_k}(\rho)$ pointwise. Therefore our choice of parameters
(c.f.$(\bigstar)$ on page \pageref{bigstar}) and the Citation
\ref{l-6.5} with the choice $N:=N_2$ and $\theta:=\theta_y$ will
guarentee that $\Gamma_{y_k}\subseteq
\Gamma_\epsilon$. Since an element fixing $y_k$ and $\epsilon$
will fix the whole ray $[y_k,\epsilon[$ pointwise, we'll also have
$\Gamma_{y_k}\subseteq \Gamma_{\varphi_\epsilon(y_k)}=\Gamma_{y_{k+1}}$.
In our argument $k$ was arbitrary, therefore 
$\Gamma_{y_{k+1}}\subseteq\Gamma_\epsilon$ and (ii).2 follows. 

It remains to proof part (iii).2: 
The vertex $x$ in question is already known to posses properties 
(i).1 and (ii).1, and the same holds true for all vertices
$\varphi_\epsilon^k(x)$ on $[x,\epsilon[$. 
$\Gamma_\epsilon$ acts transitively on
$S_\epsilon(\varphi_\epsilon^k(y))$, so we may use elements therein to
map any point in $S_1(k):=S_{\epsilon}(\varphi_\epsilon^k(y)) \cap
S_{\varphi_{\epsilon}^{k+1}(y)}(\frac{1}{2})$ 
to any other point in the same set. 
The index 
$$\left|\Gamma_{\varphi_{\epsilon}^{k+1}(y)}:
  \Gamma_{\varphi_{\epsilon}^{k+1}(y)} \cap
  \Gamma_{\varphi_{\epsilon}^{k}(y)} \right| =
\left|\Gamma_{\varphi_{\epsilon}^{k+1}(y)}:\Gamma_{\varphi_{\epsilon}^{k}(y)}
\right|$$

is therefore equal to the size of the set $S_1(k)$ which is readily
computed and finishes the proof of the lemma. \qed

Conversely, whenever we observe ever larger stabilizers, 
we approach a $\Gamma$--cuspidal end: 

\begin{lemma}\label{stab->infty=>cusp}
Suppose along some (equivalently any) ray defining an end 
of $\XG$ the order of the $\Gamma$--stabilizers of vertices 
increase to infinity. Then that end is $\Gamma$--cuspidal.
\end{lemma}
\proof
Denote the vertex sequence along the ray in question by 
$(y_i)_{i\in\mathbb{N}}$. 
The stabilizer of the edge connecting $y_i$ to $y_{i+1}$ 
equals $\Gamma_{y_i} \cap \Gamma_{y_{i+1}}$ which is of 
index at most $\q{y_{i}}+1$ in $\Gamma_{y_i}$. 
We conclude that the order of $\Gamma_{y_i} \cap \Gamma_{y_{i+1}}$ 
tends to infinity as $i$ tends to infinity as well.  

Choose an index $i_0$ such that for all $i\ge i_0$ 
$$
\#\Gamma_{y_i}>\max\{\#\Aut{B_{x_0}(\rho)}_.,\;\#\Aut{B_{x_1}(\rho)}_.\}
$$
Then the restriction map $\Gamma_{y_i}\to \Aut{B_{y_i}(\rho)}_.$ 
can not be injective, therefore using property $(\bigstar)$ we see 
that 
$$
1\neq\Gamma_{y_i}\cap\Fix{B_{y_i}(\rho)}_\Gamma.\subseteq 
\Gamma\cap gN_1g^{-1}\;.
$$
Choose an element $g_i$ in $\G k.$ mapping one of 
$x_0$, $x_1$ to $y_i$. 
According to Citation \ref{R-3.15} $\Gamma\cap g_iN_2g_i^{-1}$ 
will contain a nontrivial good unipotent element $\theta_i$.
Let $\epsilon_i$ be the end fixed by $\theta_i$. 
As alredy noticed, we may put $N:=N_2$ in Citation \ref{l-6.5}. 
With this choice of $N$ we obtain 
$\Gamma_{y_i}\subseteq\Gamma_{\epsilon_i}$. 

The parameter $n$ was chosen such that $N_2$ fixes the 
balls $B_{x_0}(1)$ and $B_{x_1}(1)$ pointwise ($(\bigstar)$). 
Hence $\theta_i\in g_iN_2g_i^{-1}$ will fix the ball $B_{y_i}(1)$, 
in particular it fixes $y_{i+1}$. 
In other words $\theta_i\in\Gamma_{y_{i+1}}$. 
As a consequence $\epsilon_i=\epsilon_{i+1}$. 
Using induction, we see that all ends $\epsilon_i$ for $i\ge i_0$ 
are equal. 
Name this end $\epsilon'$.

We claim that $\epsilon=\epsilon'$. 
Suppose this is not the case. 
Let $y_{i'}$ with $i'\ge i_0$ be the first vertex of the ray 
which lies on the line $[\epsilon,\epsilon']$. 
The group $\Gamma_{y_i}$ for $i\ge i'$ fixes $\epsilon'$ (and $y_i$), 
therfore $y_{i'}$ as well. 
This is a contradiction, since $\Gamma_{y_{i'}}$ is a finite group 
while the order of the groups $\Gamma_{y_i}$ tends to infinity. 
We conclude that $\epsilon=\epsilon'$. 
Thus $\theta_{i_0}\in\Gamma_{y_{i_0}}\subseteq\Gamma_{\epsilon}$. 
We have shown that $\epsilon$ is indeed $\Gamma$--cuspidal. \qed 

The two preceeding Lemmata together with the classification 
of groups of elliptic automorphisms of trees give us 
part~(3) of our main Theorem: 

\begin{corollar}\label{cuspidal<=>if.lf}
An end of $\XG$ is $\Gamma$--cuspidal iff its $\Gamma$--stabilizer 
is infinite and locally finite; it is then maximal with this 
property. 
Any discrete infinite and locally finite group of automorphisms 
of a locally finite tree fixes a unique end. 
In other words claim~(3) of Theorem~\ref{X mod Gamma} is true. 
\end{corollar}
\proof 
Suppose first that $\epsilon$ is an end whose 
$\Gamma$--stabilizer is infinite and locally finite.  
We know that a discrete group of tree automorphisms has 
finite point stabilizers (see \ref{tree-section}). 
Since $\Gamma_\epsilon$ is infinite, it can not have a 
common fixed point. 
It can not contain hyperbolic elements either, since 
these have infinite order. 
We conclude from Proposition 2 in \cite{CovGG} that there 
is a unique end $\epsilon'$ fixed by $\Gamma_\epsilon$, 
and that $\Gamma_\epsilon$ is the increasing union of 
the stabilizers of the vertices on any ray 
defining $\epsilon'$ and the inclusions are proper 
infinitely often. 
Uniqueness of the end gives us $\epsilon'=\epsilon$. 
This also proves that discrete infinite locally finite 
groups of automorphisms fix a unique end, hence the second 
claim is true as well.

We conclude that the order of the $\Gamma$--stabilizers 
along any ray defining $\epsilon$ tend to infinity. 
According to the previous Lemma $\epsilon$ is 
$\Gamma$--cuspidal as claimed. 
Obviously, any infinite locally finite group is contained 
in a maximal one, which fixes a unique end by what we already 
know. Therefore $\Gamma_\epsilon$ is indeed a maximal 
infinite and locally finite subgroup. 

Finally, we prove the converse of the first claim. 
Let $\epsilon$ be a $\Gamma$--cuspidal end. 
We use again Proposition 2 from \cite{CovGG}. 
This result is applicable, since $\Gamma_\epsilon$ does 
not contain any hyperbolic elements thanks to 
Corollary~\ref{SpitzenStab-}. 
We wish to prove that $\Gamma_\epsilon$ satisfies the 
third condition listed in Proposition 2 of loc.cit.; 
it will then be the (infinitely often strictly) 
increasing union of stabilizers of vertices. 
Since stabilizers of points in $\Gamma$ are finite groups, 
it will then follow that $\Gamma_\epsilon$ is infinite 
and locally finite. 

We show that $\Gamma_\epsilon$ does not satisfy the two other 
conditions listed. 
The stabilizer of an end can not contain inversions, 
therefore the second condition can not hold. 
To exclude the first possiblity it suffices to show, again 
involving discreteness, that $\Gamma_\epsilon$ is an 
infinite group. 
But this easily follows from parts (iii) and (ii) of 
Lemma~\ref{note1}. \qed

As is obvious from the last part of the previous proof 
we may reformulate again, getting 
the converse to Lemma \ref{stab->infty=>cusp}:

\begin{corollar}\label{stab->infty<=>cusp}
An end of $\XG$ is $\Gamma$--cuspidal iff the orders of the 
$\Gamma$--stabilizers of vertices along some (any) ray 
representing that end increase to infinity. 
\end{corollar}

Piecing together what we got so far, we get part~(1) and half 
of part~(2) of our main Theorem: 

\begin{corollar}\label{finish}
The quotient graph of groups defined by the $\Gamma$--action 
on $\XG'$ looks like described in part (1) of Theorem~\ref{X mod Gamma}.  
The map described in part (2) is a surjection. 
The geometric cusps correspond bijectively to the ends of any 
lift of a maximal subtree of the quotient graph. 
\end{corollar}
\proof 
Choose a fundamental $\Gamma$--transversal $Y$ with maximal subtree 
$Y_0$ in $\XG'$ (cf. \cite[I.(2.6)]{graphs} for notation). 

\emph{First step.} 
All ends of $Y_0$ are $\Gamma$--cuspidal. There are only finitely 
many of them. 

The vertices along any ray in $Y_0$ run through vertices which 
are all inequivalent modulo $\Gamma$. 
The combinatorial volume formula at the end of 
Section~\ref{rk1.geo-int} shows that the orders 
of their $\Gamma$--stabilizers must converge to infinity. 
Thanks to Lemma~\ref{stab->infty=>cusp} 
(or Corollary~\ref{stab->infty<=>cusp}) that means that the end 
defined by the ray is $\Gamma$--cuspidal. 
This proves the first statement. 
According to Proposition \ref{f of cusps} $\Gamma$ acts 
with finitely many orbits on the $\Gamma$--cuspidal ends. 
It will then map appropriate tails of rays defining ends 
lying in the same $\Gamma$--orbit 
onto each other. 
It follows that $Y_0$ does only contain finitely many inequivalent 
rays, hence has only finitely many ends. 

\emph{Second step.} 
On each ray representing an end of $Y_0$ the index condition on stabilizers 
stated in Theorem~\ref{X mod Gamma} eventually holds. 

This is immediate from the first step, using Lemma~\ref{note1}, 
part (iii). 

\emph{Third step.} 
Deleting appropriate rays representing all the ends of $Y_0$ 
leaves a finite graph. 
The geometric cusps correspond bijectively to the ends of any 
lift of a maximal subtree of the quotient graph. 

The first claim is already obvious from the first step. 
By the result of the second step, the deletion alluded to 
can be done in such a way that even the remainder of $Y$ 
is finite. 
This shows that all rays of $Y$ actually define inequivalent 
ends of the quotient graph. This observation is 
independent of $Y$, and the second claim follows. 

It remains to observe: 

\emph{Final step.} 
The map described in part (2) of Theorem~\ref{X mod Gamma} is 
a surjection. 

It may seem, that the second statement of the third step 
already takes care of part~(2) of Theorem~\ref{X mod Gamma} 
alltogether, but this is not the case. 
Lets prove the statement above. 
Given a cusp of $\Gamma$, we choose a representative 
$\Lambda$. 
It fixes a unique end $\epsilon_\Lambda$ of $\XG$, which is 
evidently $\Gamma$--cuspidal. 
Choose a ray $\mathbf{r}$ representing that end and let 
$\ol{\mathbf{r}}$ be its image in the quotient graph. 
We claim that some subray of $\mathbf{r}$ maps to a ray 
in $\Gamma\backslash\XG'$.  
Due to the structure of the quotient graph, which we have 
already established, any ray there will define an end. 
This end will be the image of the $\Gamma$--conjugacy class 
of $\Lambda$.\label{spell gp>geo} 
This spells out the definition of the map between 
cusps and geometric cusps in detail. 

The vertex stabilizers along the ray $\mathbf{r}$ have orders 
going to infinity. The stabilizers attached to any bounded 
subgraph of $\Gamma\backslash \XG'$ have bounded orders. 
Therefore the images of the vertices on the ray $\mathbf{r}$ 
must eventually leave any bounded subgraph forever. 
Thanks to the structure of the quotient graph, this means that 
these vertices converge to an end of $\Gamma\backslash \XG'$. 
Looking more closely, we see that the images of the vertices 
on a tail of the ray $\mathbf{r}$ will actually define a ray, 
thanks to the index condition established in the second step. 

The map will obviously be independent of the choice of 
$\mathbf{r}$. To see that it is independent of the choice 
of $\Lambda$, we must prove that $\Gamma$--conjugacy of 
$\Gamma_\epsilon$ and $\Gamma_{\epsilon'}$ implies that 
$\epsilon$ and $\epsilon'$ are in the same $\Gamma$--orbit. 
Both $\gamma\Gamma_\epsilon \gamma^{-1}$ and $\Gamma_{\epsilon'}$ 
are infinite, locally finite groups thanks to 
Corollary~\ref{cuspidal<=>if.lf} without fixed points. 
Involving the Proposition form \cite{CovGG} used in the 
proof of that Corollary again, we see, that the end they fix 
is unique. If they are equal, we therefore conclude that 
$\gamma.\epsilon=\epsilon'$. The map described is therefore 
well defined. 

Finally, to show surjectivity, take any end $\partial$ of 
$\Gamma\backslash \XG'$. Extend a ray representing that 
end to a maximal subtree, lift it, and extend it to a 
$\Gamma$--transversal $(Y,Y_0)$. 
The end $\epsilon$ of $Y_0$ mapping to $\partial$ is cuspidal thanks to 
the first step. One easily checks that the $\Gamma$--conjugacy 
class of the cusp subgroup of $\Gamma_\epsilon$ maps to 
$\partial$ under the map just described. We are done. 
\qed

We next attempt to prove part (2) of Theorem~\ref{X mod Gamma}. 
It amounts to showing that, whenever two vertex sequences 
$(x_i)$ and $(y_i)$ define $\Gamma$--cuspidal rays, 
such that there is a sequence of elements $(\gamma_i)$ in $\Gamma$ 
with $\gamma_i.x_i=y_i$, there is actually a single element 
$\gamma$ in $\Gamma$ mapping the corresponding ends on each other. 

We will show that this is the case using a concept, 
which allows the geometric interpretation of the contraction 
of cusp construction we promised earlier. 
We are going to build a quotient space of the tree $\XG'$ 
by identifying horoballs to points. These will correspond 
to the $\Gamma$--orbits of the contracted rays. 
When forming a quotient of a space acted upon by a group, 
one may lose control over stabilizers, unless further 
conditions are imposed. 
One will be able to avoid this problem if the following 
conditions are met: 

\begin{definition}\label{precisely-invariant}
Let $\Gamma$ be a group acting on the set $X$. 
A subset $Y$ of $X$ is called precisely invariant under 
the subgroup $\Lambda$ of $\Gamma$ iff the 
conditions 
\begin{itemize}
\item[(1)] $\Lambda=\Stab Y_\Gamma.$ and 
\item[(2)] $g(Y)\cap Y=\emptyset$ for all $g\in \Gamma\setminus \Lambda$ 
           hold. 
\end{itemize}
An $n$--tuple $(Y_1,\ldots,Y_n)$ of subsets of $X$ is called 
precisely invariant under an $n$--tuple 
$(\Lambda_1,\ldots,\Lambda_n)$ of subgroups iff 
\begin{itemize}
\item[(1)] $Y_i$ is precisely invariant under $\Lambda_i$ for all $i$ and 
\item[(2)] for all $g\in\Gamma$ and indices $m\neq k$ we have 
           $g(Y_m)\cap Y_k=\emptyset$. 
\end{itemize}
\end{definition}

Borrowing a notion from hyperbolic geometry, 
we will call a horoball $B_\epsilon(x)$ in a tree independent 
(with respect to a group of automorphisms $\Gamma$) iff 
it is precisely invariant under $\Gamma_\epsilon$ in $\Gamma$. 

The following result extends Lemma 2.5 from \cite{l.SL2} 
to all semisimple group of rank 1. 
Via Lemma~\ref{note1} this result depends again crucially 
on Citation~\ref{l-6.5}. 

\begin{lemma}[Existence of independent horoballs] 
Let $\Gamma$ be a lattice in $\G k.$, where $\G.$ is connected semisimple $k$--group of $k$--rank 1. 
\begin{itemize}
\item[(i)] Let $x_\epsilon$ be a vertex of $\XG'$, such that any 
        vertex $x$ in $B_\epsilon(x_\epsilon)$ contains 
        a nontrivial good unipotent $\theta_x\in \Gamma_\epsilon$ 
        fixing $B_x(\rho)$ pointwise (; existence of $x_\epsilon$ 
        is assured by part (iii).1 of Lemma \ref{note1}). 
        Then $B_\epsilon(x_\epsilon)$ is independent with respect 
        to $\Gamma$. 
\item[(ii)] If $\epsilon$ and $\epsilon'$ are $\Gamma$--inequivalent 
            $\Gamma$--cuspidal ends, then 
            $(B_\epsilon(x_\epsilon),B_{\epsilon'}(x_{\epsilon'}))$ 
            is precisely invariant under 
            $(\Gamma_\epsilon,\Gamma_{\epsilon'})$. 
\end{itemize}
\end{lemma}
\proof 
We defer the proof of (i). Modulo (i) claim (ii) and the 
verification of the property (2) for claim (i) come down to 
the following: 
Given $y\in B_\epsilon(x_\epsilon)$ and $y'\in B_{\epsilon'}(x_{\epsilon'})$ 
and further an element $\gamma\in\Gamma$ with $\gamma.y=y'$, we will 
have $\gamma.\epsilon=\epsilon'$. 
To see this, choose a nontrivial good unipotent $\theta_y\in\Gamma_y$, 
whose conjugate under $\gamma$ lies in 
$\Gamma_{y'}\subseteq\Gamma_{\epsilon'}$. 
This is seen to be possible by applying part (iii).1 of Lemma \ref{note1}
twice. 
We have $\gamma\theta_y\gamma^{-1}.\epsilon'=\epsilon'$. 
Since a nontrivial good unipotent element fixes a unique end, 
we infer $\gamma.\epsilon=\epsilon'$, showing the claim. 

We still have to prove that property (1) holds in claim (i). 
It suffices to show that $\Gamma_\epsilon$ leaves the horoball 
$B_\epsilon(x_\epsilon)$ invariant. This is the case, since 
we know from Corollary~\ref{SpitzenStab-}, that all elements of 
$\Gamma_\epsilon$ are elliptic.
\qed

The whole point of this notion is that it allows the following 
result to be proved: 

\begin{lemma}\label{geo-cusp.contr}
Let $\Gamma$ be a nonuniform lattice in $\G k.$, where 
$\G.$ is connected semisimple $k$--group of $k$--rank 1 
over the local field $k$. 
Let $\{\epsilon_1,\ldots, \epsilon_c\}$ be a system of pairwise 
inequivalent representatives of the $\Gamma$--orbits of 
$\Gamma$--cuspidal ends. 
For each $\epsilon_i$ choose an independent horoball $B_i$ 
around $\epsilon_i$ in $\XG'$. 
Identify each of the $\Gamma$--translates 
of any $B_i$ to a point, i.e., form the quotient graph 
$\ol{X}_{\G.}:= \XG'/\Gamma B_1\cup \cdots \cup \Gamma B_c$ 
(for notation cf. \cite[p.80]{trees}). 
It is a tree thanks to Proposition 13 from \cite[p.22]{trees}.  
We call $\ol{X}_{\G.}$ a tree obtained by contracting horoballs 
around cusps. 

Denote the vertex, which is the image of $B_i$ in 
$\ol{X}_{\G.}$ by $\ol{B}_i$. 
We have $\Gamma_{\ol{B}_i}=\Gamma_{\epsilon_i}$. 
The contraction map $\XG'\to \ol{X}_{\G.}$ induces a map 
$\Gamma\bbackslash \XG'\to \Gamma\bbackslash \ol{X}_{\G.}$, 
which is just the contraction of cusps construction described 
earlier. The tail of a ray $r(i)$ (representing a geometric cusps of 
$\Gamma$) which gets contracted to $\Gamma \ol{B}_i$ is determined by 
the condition that some ($\iff$ any) lift of its initial point lies 
on the boundary horosphere of a $\Gamma$--translate of $B_i$. 

Since the points  $\Gamma \ol{B}_i$ are different for different 
indices, we get as a Corollary that the canonical map from the cusps 
of $\Gamma$ to the geometric cusps of $\Gamma$ is injective. 
\end{lemma}
\proof 
Digesting the definitions and constructions involved, the only nontrivial 
claim to prove is disjointness of horoballs of the form 
$yB_i$ and $y'B_{i'}$ whenever $i\neq i'$ or $i=i'$ but 
$y^{-1}y'\notin \Gamma_{\epsilon_i}=\Gamma_{\epsilon_{i'}}$. 
But this has been proved in the above Lemma. \qed

